\documentclass[a4paper,10pt,reqno]{amsart}

\input{psfig.tex}   
\input{psfig.sty}
   
\usepackage{latexsym}
\usepackage{amssymb} 
\usepackage{amsbsy} 
\usepackage{amsthm}
\usepackage{amsfonts}   
\usepackage{amscd}
\usepackage{amsmath}     
\usepackage{yfonts}

\def\Z{{\mathbb Z}}   
\def\C{{\mathbb C}}   

\newcommand{\bu}{b_{1}}   
\newcommand{\bd}{b_{2}}   
\newcommand{\bt}{b_{3}}   
\newcommand{\bj}{b_{j}}   
\newcommand{\bjmu}{b_{j-1}}   
\newcommand{\bc}{b_{1}}   
\newcommand{\bm}{b_{2}}   
\newcommand{\be}{\beta}   
\newcommand{\al}{\alpha}

\newcommand{\mct}{\mathcal{T}}      
\newcommand{\bn}{B_{n}}   
\newcommand{\hqn}{H(Q,n)}   
\newcommand{\vsp}{\vspace{10pt}}   
\theoremstyle{plain}   
\newtheorem{theorem}{Theorem}[section]   
\newtheorem{proposition} {Proposition}[section]   
\newtheorem{corollary} {Corollary}[section]   
\newtheorem{lemma}{Lemma}[section]   
\newtheorem{conjecture}{Conjecture}[section]

\theoremstyle{remark}
\newtheorem{example}{Example}[section]
\newtheorem{remark}{Remark}[section]
\theoremstyle{definition} 
\newtheorem{definition} {Definition}[section] 

\begin{document}

\title[Polynomial invariants of links]{Polynomial invariants of links satisfying cubic 
skein relations}

\author[P.Bellingeri]{Paolo Bellingeri}
\address{Institut Fourier, BP 74, Univ.Grenoble I, 
Math\'ematiques, 38402 Saint-Martin-d'H\`eres cedex, France}
\curraddr{D\'epartement de Math\'ematiques,
BP 051, Univ. Montpellier II, 34095 Montpellier cedex, France}
 
\email{bellinge@fourier.ujf-grenoble.fr, 
bellinge@math.univ-montp2.fr}
 \urladdr{http://www-fourier.ujf-grenoble.fr/$\sim$bellinge}

\author[L.Funar]{Louis Funar}
\thanks{Partially supported by a Canon grant}
\address{Institut Fourier, 
BP 74, Univ.Grenoble I, Math\'ematiques,  
38402 Saint-Martin-d'H\`eres cedex, France} 
\email{funar@fourier.ujf-grenoble.fr}
\urladdr{http://www-fourier.ujf-grenoble.fr/$\sim$funar}
    
\date{}

\begin{abstract}   
  The aim of this paper is to define two link invariants    
satisfying cubic skein relations.   
In the hierarchy of polynomial invariants determined  by explicit 
skein relations they  are the next level  of complexity   
after Jones, HOMFLY, Kauffman and Kuperberg's $G_2$ quantum   
invariants.    
Our method consists in the study of Markov traces on a suitable tower   
of quotients    
of cubic Hecke algebras extending Jones approach.   
 \end{abstract}   
 

\maketitle

\section{Introduction}

\subsection{Preliminaries}   
   
J.Conway showed that the Alexander polynomial of a knot,   
when suitably normalized, satisfies the following skein relation:

$$  \nabla \left(   \raisebox{-4mm}{\psfig{figure=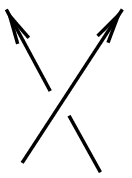,width=20pt}}   \right)   
 - \nabla \left(   \raisebox{-4mm}{\psfig{figure=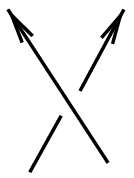,width=20pt}}   \right) =  
(t^{-1/2} - t^{1/2}) \nabla \left(   \raisebox{-4mm}{\psfig{figure=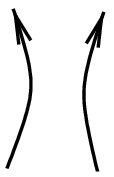,width=20pt}}   \right)  
$$

Given a knot diagram, one can always change some of 
its crossings    
such that the modified diagram represents the unknot.    
Therefore, one can use the skein relation for a recursive 
computation    
of $\nabla$, although this algorithm is rather time consuming,     
since it is exponential.

In the mid eighties V.Jones discovered another invariant verifying a different    
but quite similar skein relation, namely:

$$  t^{-1} V \left(   \raisebox{-4mm}{\psfig{figure=des.eps,width=20pt}}   \right)   
 - t \, V \left(   \raisebox{-4mm}{\psfig{figure=sin.eps,width=20pt}}  
 \right)   
=  
(t^{-1/2} - t^{1/2}) V \left(   \raisebox{-4mm}{\psfig{figure=inf.eps,width=20pt}}   \right)
$$

\noindent  which was further generalized to a 2-variable invariant by replacing   
the factor $(t^{1/2}-t^{-1/2})$ with a new variable $x$. The latter  
one was shown to specialize to both Alexander and Jones 
polynomials.

The Kauffman polynomial is another extension of Jones   
polynomial which satisfies     
a skein relation, but this time 
in the realm of unoriented diagrams.    
Specifically, the formulas:

$$ \Lambda \left(   \raisebox{-4mm}{\psfig{figure=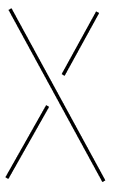,width=17pt}}  
  \right)   
+ \Lambda \left(   \raisebox{-4mm}{\psfig{figure=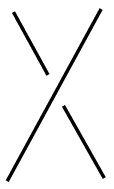,width=17pt}}  
  \right) =  
z  \left( \Lambda \left(     
\raisebox{-4mm}{\psfig{figure=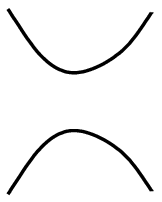,width=20pt}}   \right)  
+ \Lambda \left(   \raisebox{-4mm}{\psfig{figure=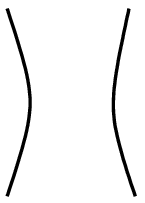,width=20pt}}  
  \right) \right)  
$$  
  
$$ \Lambda \left(  \raisebox{-2mm}{\psfig{figure=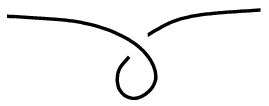,width=50pt}}  
  \right)  =  
a \Lambda \left(   \raisebox{0.2mm}{\psfig{figure=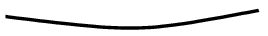,width=50pt}}  
  \right)  
$$

\noindent  define a regular isotopy invariant of links which 
can be  renormalized, by using the writhe of the oriented diagram, 
in order to become a link invariant.    
Remark that some elementary manipulations show that 
$\Lambda$ verifies the following  cubical skein relation:

$$ \Lambda \left( 
  \raisebox{-5mm}{\psfig{figure=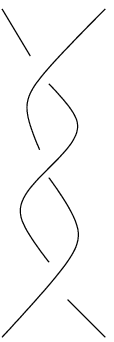,width=14pt}}    
\right) =  
\left( \frac{1}{a} + z \right)\,  \Lambda \left(   \raisebox{-4.5mm} 
{\psfig{figure=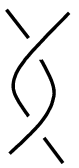,width=16pt}}   \right) -  
\left( \frac{z}{a} + 1 \right) \,   \Lambda \left(   \raisebox{-4mm}{\psfig{figure=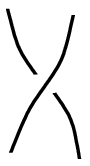,width=17pt}}   \right)  
+ \left( \frac{1}{a}\right) \, \Lambda \left(   \raisebox{-4mm}{\psfig{figure=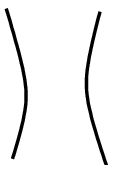,width=21pt}}   \right)  
$$

It has been recently proved 
(\cite{DP}, and Problem 1.59 \cite{kir}) that  this   
relation alone is not  sufficient for a   
recursive computation  of $\Lambda$. Whenever the skein  
relations and the value of the invariant for the unknot 
are sufficient to determine its values for all links,
the system of skein  relations will be said to be complete.   
Several results concerning the incompleteness of 
higher degree unoriented skein relations and their skein modules 
have been obtained by J.Przytycki and his 
students (see e.g. \cite{DP,Prz2,Prz3}). 
   
These invariants were generalized to quantum invariants 
associated to general Lie  algebras, super-algebras 
and their representations. V.Turaev   
(\cite{Turaev}) identified the HOMFLY and Kauffman 
polynomials with  the invariants obtained from the 
series $A_n$ and $B_n, C_n, D_n$   
respectively. G.Kuperberg (\cite{Ku}) defined the $G_2$-quantum   
invariant of knots  by means of skein relations, by  
making use of trivalent    
graphs diagrams and exploited further 
these ideas in \cite{Ku2}, for  spiders of   rank 2 Lie algebras. 
The skein relations satisfied by the    
quantum invariants coming from simple Lie algebras 
were approached  also via weight systems and the Kontsevich 
integral in  \cite{LM1,LM2}, for the classical series, 
and in  \cite{BS1,BS2} for the case of
the Lie algebra $\textgoth{g}_2$ of $G_2$.

Notice that any link invariant coming from 
some R-matrix $R$ verifies a skein relation of the type:

$$ \sum_{j=0}^{n}  a_j \left\langle \; \left. \raisebox{-9mm}{\psfig{figure=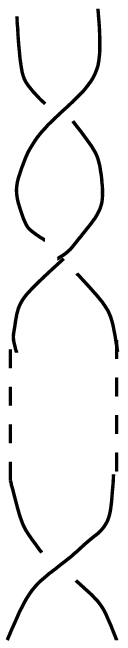,width=10pt}}   
 \right\}  
\; \mbox{j twists} \;  
\right\rangle =0  
$$

\noindent  which can be derived from the polynomial equation 
satisfied by the R-matrix $R$.     
    
Let us mention that the skein relations are somewhat related to    
the representation theory of the Hopf algebra associated to 
the R-matrix $R$.  In particular, there are no other known 
invariants given by means of a complete family of 
skein relations, but those from above. 
Moreover, one expects that the quantum invariants 
associated to other Lie (super) algebras or by cabling 
the previous ones should satisfy skein relations of 
degree at least 4, as already the
Kuperberg $G_2$-invariant does.

This makes the search for an explicit set of complete skein 
relations,  in which at least one relation is cubical, 
particularly difficult and  interesting. 
This problem was first considered in \cite{Fun} and   
solved in a particular case. The    
aim of this paper is to complete the result of \cite{Fun} by    
constructing a deformation of the previously constructed    
quotients of the cubic Hecke algebras  
and of the  Markov traces supported   
by these algebras. We obtain in  this way two link  invariants, 
denoted by $I_{(\al,\,\be)}$   
and  $I^{(z,\, \delta)}$, which are  
recursively computable and uniquely determined by 
two skein relations. 
Explicit computations show that $I_{(\al,\,\be)}$  
detects the chirality of the knots 
with  number crossing   at most 10  
where HOMFLY,  Kauffman and their $2$-cablings fail. 
On the other hand, as HOMFLY,  Kauffman and  
their $2$-cablings (\cite{LL,Prz}), it seems that 
our invariants do not  distinguish between mutant knots. 
We recall that the some mutant knots can be 
distinguished by the $3$-cablings of the HOMFLY polynomial 
(see \cite{Mur}).

\vspace{0.5cm}  
   
\noindent  {\bf Acknowledgements.} Part of this work was done during the second author's visit to the    
Tokyo Institute of Technology, whose support and hospitality are   
gratefully acknowledged. The authors are thankful to      
Christian Blanchet, Emmanuel Ferrand, Thomas Fiedler, 
Louis Kauffman, Teruaki Kitano,  Sofia Lambropoulou, 
Ivan Marin, Jean Michel, Hugh R. Morton, Luis Paris and  
Vlad Sergiescu for useful discussions, remarks and  suggestions.

\subsection{The main result}

The aim of this paper is to define two  link invariants    
by means of a complete set of skein relations. 
More precisely we will prove   
the following Theorem (see section 5):

\begin{theorem} There exist a link invariant  
$I_{(\al,\,\be)}$  which is uniquely determined by the   
two skein relations shown in (1) and (2) and 
its value for the unknot, which is traditionally 1.    

\vspace{-0.5cm}   
\begin{eqnarray}
 \raisebox{-7mm}{\psfig{figure=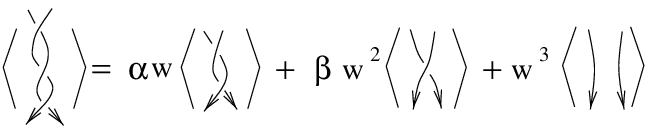,width=8cm}}
\end{eqnarray}

\begin{eqnarray}
 \raisebox{-70mm}{\psfig{figure=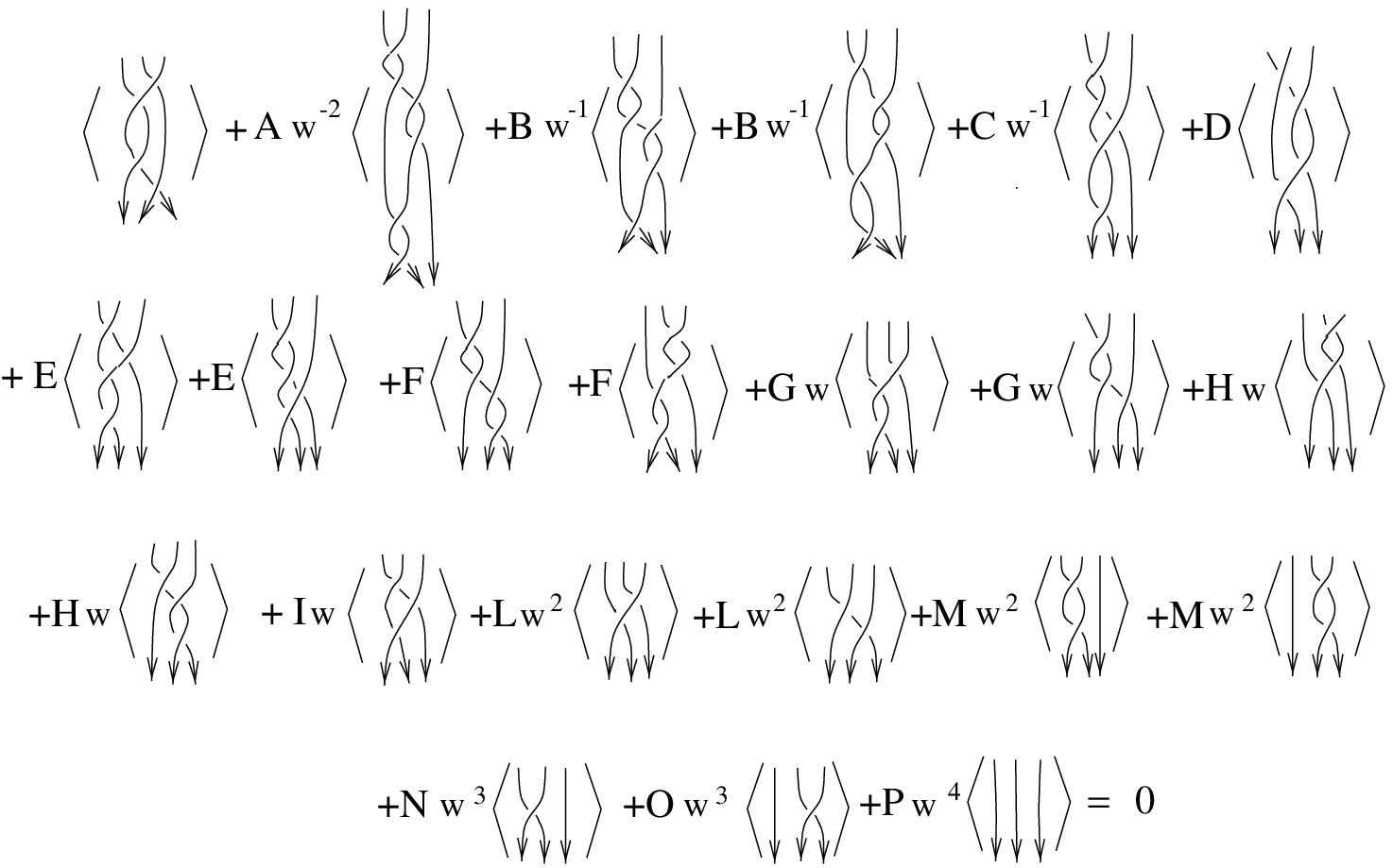,width=13cm}}
\end{eqnarray}

\noindent The invariant takes values in:    
   
$$   
\frac{\Z [\al, \, \be , \,(\be^2 - 2\al)^{\pm \epsilon/2} , \, (\al^2 + 2\be)^{\pm \epsilon/2} ]}   
{({H}_{(\al,\,\be)})}      
$$   
\noindent  where  $\epsilon-1\in\{0,1\}$ is  
the number of link components modulo 2, and  
${H}_{(\al,\,\be)}$ is the following polynomial: 
 
\[{{H}_{(\al,\,\be)}} =  8\alpha^6- 8\alpha^5\beta^2+2\alpha^4    
   \beta^4+36\alpha^4\beta -34\alpha^3\beta^3+17\alpha^3+8\alpha^2     
 \beta^5+32\alpha^2\beta^2 -36\alpha\beta^4    
 +38\alpha\beta+8\beta^6-17\beta^3+8\]    
Here $(Q)$ denotes the ideal generated by the element $Q$ in the algebra   
under consideration.
The values of the polynomials $A, B, C,\ldots, P$ 
appearing in the skein relations for  $I_{(\al,\,\be)}$ are  
given in the table below:

$$   
\begin{array}{|l|l|}   
\hline   
 w=((\alpha^2+2\be)/(2\alpha-\be^2))^{1/2} & A= (\be^2 - \al) \\   
\hline    
 B=(\al^2 - \al \be^2 - \be) & C=(\al^2 - \al \be^2)  \\   
\hline   
D=(1+ 2\al \be  + \al^2 \be^2 - \al^3)  & E=(1+ \al \be  + \al^2 \be^2 - \al^3)\\   
\hline   
F=(1+ 2\al \be  - \be^3)  & G=(\al \be^3 - 2 \al - 2 \al^2 \be) \\   
\hline   
H= (\al \be^3 - 2 \al - 2 \al^2 \be + \be^2) & I=(\al^4 -\al^3 \be^2 - 2\al^2\be - 3 \al) \\   
\hline   
 L=(2 \al^3 \be + 3\al^2 - \al^2 \be^3 - \al \be^2) & M=   
(\be^4-2\be-3\al\be^2+\al^2)  \\   
\hline   
 N= (1+ 4\al \be +3\al^2 \be^2 - \al^3 - \al \be^4- \be^3) &   
O=(1+ 3\al \be + 3\al^2 \be^2 - \al^3 - \al \be^4) \\   
\hline   
   P=( 3\be^2 - \be^5 - 2 \al - 3 \al^2 \be + 4 \al \be^3 ) &   
   \\   
\hline      
\end{array}  
$$      
\begin{center}  
{Table 1}  
\end{center}    
  
Furthermore there exists a second link invariant 
$I^{(z,\delta)}$, which is determined by the skein relations 
(1) and (2), but with another set of coefficients. Specifically, 
$I^{(z,\delta)}$ takes values in:  
   
$$ \frac{\Z [z^{\pm \epsilon/2}, \delta^{\pm \epsilon/2 }]}   
{({P}^{(z,\, \delta)})}      
$$   
   
\noindent  where  $\epsilon$ is as above and   
${P}^{(z,\, \delta)}$ is 
the following polynomial:     
\[ {P^{(z,\, \delta)}}=z^{23}+z^{18} \delta-2 z^{16}  
\delta^2-z^{14} \delta^3-2 z^9 \delta^4+2 z^7 \delta^5+\delta^6  
z^5+\delta^7\]   
The values of the rational functions $A, B, C,\ldots, P$ 
corresponding to the skein relations 
of  $I^{(z,\delta)}$ are obtained    
as follows:  set first  
$w=\left(-\frac{z^3}{\delta}\right)^{1/2}$   
and  substitute further 
$\al =-\frac{z^7+\delta^2}{z^4 \delta}$ and    
$\be= \frac{\delta-z^2}{z^3}$ in the other entries of  table 1.    
\end{theorem}

\subsection{Properties of the invariants}   
 
The following summarize the main features of these 
invariants (see section 6):  
\begin{enumerate}   
\item  they distinguish all  knots  with number crossing   
 at most 10 that have the same HOMFLY   
polynomial, and thus they are independent from HOMFLY.   
However, like HOMFLY and Kauffman polynomials, they seem to not   
distinguish among mutants knots. 
In fact, they do not distinguish between the  Kinoshita-Terasaka knot
and the Conway knot, which are the simplest non-equivalent mutant 
knots.   
\item $I_{(\al,\,\be)}=I_{(-\be,\,-\al)}$  for amphicheiral knots, and    
$I_{(\al,\,\be)}$ detects the chirality of all those 
knots with number  crossing at most 10, whose 
HOMFLY, Kauffman polynomials as well as the 
$2$-cabling of HOMFLY fail to detect.   
\item $I_{(\al,\,\be)}$ and  $   
I^{(z,\, \delta)}$ have a {\em cubical} behaviour.   
\end{enumerate}   
Let us explain briefly what we meant by {\em cubical behaviour}.    
\begin{definition}    
A Laurent polynomial $\sum_{j \in \Z} c_{j}a^{j}$ is a   
$(n,k)$-polynomial (for $n,k \in \Z$)  
if  $c_j = 0$ for $j \ne k \,({\rm modulo }\; n)$.   
\end{definition}

\begin{remark}   
\begin{enumerate}  
\item   
The HOMFLY polynomial $V(l,m)$ can be written as    
$\sum_{k \in\Z} R_k(l)m^k$ and respectively as   
$\sum_{k \in\Z} S_k(m)l^k$, where $R_k(l)$ and   
$S_k(m)$  are $(2,k)$-Laurent polynomials   
fulfilling $R_{2k+1}(l)=S_{2k+1}(m)=0$.   
\item The  Kauffman polynomial   
can be written as $\sum_{k \in\Z} U_k(l)m^k$  and 
respectively as   
$\sum_{k \in \Z} T_k(m)l^k$, where $U_k(l)$ and $T_k(m)$    
are $(2,k+1)$-Laurent polynomials.   
\end{enumerate}  
\end{remark}  
In this respect the HOMFLY and Kauffman polynomials have a quadratic   
behaviour.    
\begin{proposition}   
$I_{(\al,\,\be)}$ and    
$I^{(z,\, \delta)}$ have a cubical behaviour, i.e. for each link $L$   
there exists some $l\in\{0,1,2\}$   
so that:
$$I_{(\al,\,\be)}(L)= \frac{\sum_{k \in \Z_+}    
P_{k}(\be)\al^{k}}{\sum_{k \in \Z_+} Q_{k}(\be)\al^{k}}=   
\frac{\sum_{k \in \Z_+} M_{k}(\al)\be^{k}}{\sum_{k \in \Z_+}   
N_{k}(\al)\be^{k}}$$    
where $P_{k}, Q_k, M_{k}, N_{k}$ are  $(3, k+l)$-polynomials,  and    
 $$   
I^{(z,\, \delta)}(L)= \sum_{k \in \Z} H_{k}(\delta)z^{k}= \sum_{k \in \Z} G_{k}(z)\delta^{k}$$   
 where $H_{k}, G_{k}$ are  $(3, k)$-Laurent polynomials.   
\end{proposition}

\subsection{Comments}  
There are three link invariants  coming   
from Markov traces on cubic Hecke algebras, presently known. 
First, for each quadratic factor $P_i$ of the cubic 
polynomial $Q$ one has a Markov  trace which factors 
through the usual Hecke algebra $H(P_i, n)$, 
yielding a re-parameterized  HOMFLY invariant. 
Then there is the  Kauffman polynomial   
and the invariant $I_{(\alpha,\beta)}$ (or $I^{(z,\delta)}$)
introduced in the present paper.   
It would be very interesting to find whether there exists 
some relationship between them. The explicit computations below 
show that the new invariants are independent on HOMFLY, 
Kauffman and their 2-cablings. 

Further, one expects that our invariants belong to a family  
of genuine two-parameter invariants, as expressed in the 
following:   
\begin{conjecture}   
There exists a Markov trace on $H(Q,n)$ taking values in an 
algebraic  extension of $\Z [\al, \, \be]$, which 
lifts the Markov trace  underlying $I_{(\alpha,\beta)}$.     
\end{conjecture}  
In other words, the non-determinacy $H_{(\al,\be)}$ in  
$I_{(\alpha,\beta)}$  can be removed. Notice that the polynomials 
$H_{(\al,\be)}$ and $P^{(z,\delta)}$  define irreducible planar   
algebraic curves which are not rational. In particular, 
one cannot express  explicitly the invariants as one 
variable polynomials.

\subsection{Cubic Hecke algebras}

The form of the first skein relation (1) 
explains the appearance of  cubic quotients 
of braid group algebras 
$\C [ B_n ]$. 
Recall that the braid group $B_n$ on $n$ strands
is given by the presentation: 
\[B_n= \langle b_1,\dots, b_{n-1} \,  |  \,  b_i b_j = b_j b_i , \, |i-j|>1 \; \; \mbox{and} \; \;   b_i b_{i+1} b_i= b_{i+1} b_i b_{i+1},   
\, i<n-1 \, \rangle\] 
Furthermore we define the  cubic   
Hecke algebra by analogy with the usual (i.e. quadratic) 
Hecke algebra (see \cite{Bou}), as follows:   
\[   
\hqn= \C[\bn] / (Q( \bj)\,; \; j=1, \dots, n-1) 
\]   
   
\noindent  
where 
$Q(\bj)=\bj^3- \alpha \, \bj^2 - \be \,\bj -1$  
is a cubic polynomial, 
which will be fixed through out this paper.

Our purpose is to construct  Markov traces on the 
tower of cubic  Hecke algebras since these will eventually lead 
to link invariants.  This method was pioneered by V.Jones 
(\cite{Jones}) and A.Ocneanu, who applied it to 
the case of usual Hecke algebras and obtained the celebrated 
HOMFLY polynomial. Later on several authors 
(see \cite{GL,I,La,O}) employed more sophisticated  
algebraic and combinatorial tools in searching 
for Markov traces on other Iwahori-Hecke algebras, 
for instance those of type ${B}$  which are
leading to invariants for links in a solid torus.

The cubic Hecke algebras are particular cases of the generic 
cyclotomic Hecke algebras, introduced by M.Brou\'e and G.Malle  
(see \cite{BMM1}) and studied in \cite{BMM2,BMM}, in connection 
with  braid group representations.   
Recall the following  
results concerning the structure of the 
cyclotomic Hecke algebras with $Q(0)\neq 0$ (according to 
\cite{BMM1,BMM2,BMM,Cox} and  \cite{Cox2}, p.148-149):    
   
\begin{enumerate}    
\item $\dim_{\C}H(Q,3)=24$ and $H(Q,3)$ is isomorphic to  the   
 group algebra of the binary tetrahedral group 
$\langle2,3,3\rangle$    
of order 24, i.e. the linear group $SL(2, {\Z}/3\Z)$.    
\item $\dim_{\C}H(Q,4)=648$ and $H(Q,4)$ is the   
group algebra of  the finite group 
$G_{25}$ from the Shepard-Todd classification (see  
\cite{she}).    
\item $H(Q,5)$ is the cyclotomic Hecke algebra of the  
group $G_{32}$, whose order is $155520$. It is conjectured that this 
algebra is free of finite dimension which would imply 
(by using the Tits 
deformation theorem) that it is isomorphic to the group algebra  
of $G_{32}$.  
\item $\dim_{\C}H(Q,n)=\infty$ for $n\ge 6$.   
\end{enumerate}   
   
Thus a direct definition of the trace on $H(Q,n)$ for 
$n\geq 6$ is  highly a nontrivial matter, because  
it would involve in particular, the  explicit solution of 
the conjugacy problem in these algebras.    
In order to circumvent these difficulties
one introduces a tower of smaller quotients  
$K_n (\al,\be)$ by adding  one more relation to $H(Q,3)$,
as follows:    
   
\[ b_2b_1^2b_2+ R_0=0 \]
where 
\begin{eqnarray*} 
R_0 & = &   
 A \,\, \bc^2 \, \bm^2 \, \bc^2  +     
 B\,\, \bc \, \bm^2  \, \bc^2 +   
B \,\, \bc^2 \, \bm^2  \, \bc        
+ C\,\, \bc^2 \, \bm \, \bc^2 +   
D \,\, \bc   \, \bm^2  \, \, \bc +     
E\,\, \bc   \, \bm   \, \bc^2  +     
E\,\, \bc^2    \, \bm   \, \bc +  \\ 
& & 
F\,\, \bm^2   \, \bc^2  +     
 F\,\, \bc^2   \, \bm^2  +     
G \,\, \bm   \, \bc^2 + 
G\,\, \bc^2   \, \bm +   
H\,\, \bm^2   \, \bc  +     
H\,\, \bc   \, \bm^2 +   
I\,\, \bc \, \bm \, \bc +   
L \,\,\bm \, \bc +   \\ 
 & & 
L \,\, \bc \, \bm +   
M \,\,\bc^2 +   
M \,\, \bm^2  +   
N \,\,\bc +   
O \,\,\bm +   
P
\end{eqnarray*}
and $A, B, \dots, P$ are the functions from table 1.   

\begin{remark}  
The main feature of these quotients is the fact that 
the algebras $K_n(\al,\be)$ are finite  
dimensional for all values of $n$. Moreover, these 
algebras do not collapse for large $n$, thus yielding 
an interesting tower of algebras.  
\end{remark}   
   
\begin{remark}  
Let us explain the heuristics behind that choice for 
the additional    relation. For generic $Q$ the algebra 
$H(Q,3)$ is semi-simple   
and decomposes as   $\C^3 \oplus M_2^{\oplus3} \oplus M_3$,   
where  $M_m$ is the algebra of $m \times m$ matrices.   
The quadratic Hecke algebra $H_q(3)=\C[B_3]/(b_i^2+(1-q)b_i-q)$  
arises as a quotient of $H(Q,3)$ by killing 
the factor $\C \oplus M_2^{\oplus2} \oplus M_3$. 
It is known that Jones and HOMFLY polynomials can be   
derived from the unique Markov trace on the homogeneous tower 
$\cup_{n=1}^{\infty}H_q(n)$.   
In a similar way,  the rank 3 Birman-Wenzl algebra (\cite{BW}) 
- which supports an unique Markov trace 
inducing the Kauffman polynomial - is the quotient of $H(Q,3)$  
by the factor $\C\oplus M_2^2$.   
In our case we introduced  the extra relation above
which kills precisely the central factor  $\C^3$ of $H(Q,3)$.     
\end{remark}  
The geometric interpretation of these relations is now obvious:   
the first skein relation (1) is the cubical relation  
corresponding to taking the   
quotient $H(Q,n)$ while  the main skein relation (2) 
defines the smaller quotient algebras  $K_n (\al,\be)$.

Our main theorem  is a consequence of the more technical result    
below (see sections 2, 3 and 4).

\begin{theorem} There are precisely four values of    
$(z, \, \bar{z})$ (formal expressions in $\al$ and $\be$) 
for which there exists a
Markov trace $\mct$  on  the tower $\cup_{n=1}^{\infty}K_n (\al,\be)$ 
with parameters $(z, \, \bar{z})$ i.e. verifying
the following conditions:   
\begin{enumerate}   
\item  $\mct (xy) = \mct (yx) \mbox{ for all } x, y\in 
K_n(\al,\be) \mbox{ and all } n$,    
\item  $\mct (x b_{n-1}) = z \mct (x)  \mbox{ for all } x\in 
K_n(\al,\be) \mbox{ and all } n$,   
\item  $\mct (x b_{n-1}^{- 1}) = \bar{z} \mct (x) 
 \mbox{ for all } x\in 
K_n(\al,\be) \mbox{ and all } n$.   
\end{enumerate}    
The first pair $(z,  \,\bar{z})$ is given by:   
\[ z=\frac{2\alpha-\be^2}{\alpha \be + 4}, \,  
\bar{z}=-\frac{\alpha^2+2\be}{\alpha \be + 4} \]   
and the corresponding trace takes values as follows:  
\[\mct_{\al, \, \be } : \cup_{n=1}^{\infty}K_n (\al,\be) \to 
\frac{\Z [\al, \, \be , \, (\al \be + 4)^{-1}]}{   
(H_{(\al,\,\be)})}\]   
The other three solutions are not rational functions on 
 $\al$ and $\beta$, but nevertheless one can express 
$\al, \, \beta$ 
and $\bar{z}$  as rational functions   
of  $z$ and $\delta$, where $\delta= z^2 (\beta z + 1)$. 
Specifically, we have a Markov trace:
\[ \mct^{(z,\, \delta)}: \cup_{n=1}^{\infty}K_n(\al, \, \be)  \to 
\frac{\Z [z^{\pm 1}, \delta^{\pm 1}]}{(P^{(z,\, \delta)})} \]  
where:
\[ \be= \frac{\delta-z^2}{z^3}, \;   \al=-\frac{z^7+\delta^2}
{z^4 \delta} \;  
\mbox{ and } \bar{z}=-\frac{z^4}{\delta}\]   
\end{theorem}   
   
\begin{remark}
For particular values of $(\alpha,\be)\in \C$ one might find 
that the indeterminacy ideal for the respective Markov traces 
is smaller than the specialization of the ideal above. 
A specific example is the $\Z/6\Z$-valued invariant, 
corresponding 
to the values $\alpha=\beta=0$ in \cite{Fun}, which is a
specialization of the invariant $I^{z,\delta}$ for
$z^3=-1$ and $\delta=z^2$. We can refine the
general Markov trace in order to restrict to a  $\Z/3\Z$-valued
trace (see section 6), but this refinement does not 
survive the deformation process.   
\end{remark}

There is a natural way to convert a Markov trace 
$\mct$ into a link  invariant, by setting:    
   
\[I(x)= \left (\frac{1}{z\bar{z}}\right )^{\frac{n-1}{2}}\left(\frac{\bar{z}}{z}\right)^{\frac{e(x)}{2}}\mct(x)\]   
where $x\in B_n$ is a braid representative of the link $L$ and   
$e(x)$ is the exponent sum of $x$.
    
Therefore we derive two invariants $I_{(\al,\,\be)}$ and  $   
I^{(z,\, \delta)}$ from the previous Markov traces, which 
satisfy the claimed skein relations.

\subsection{Outline of the proof}

We will prove by recurrence on $n$ that a Markov trace on   
$K_{n}(\al,\be)$   
extends to a Markov trace on $K_{n+1}(\al,\be)$.    
Since there is a nice system of generators for $K_{n+1}(\al,\be)$  
constructed inductively  starting from a 
generators system for $K_{n}(\al,\be)$, 
such an extension, whenever it exists, it must be unique.   
This is a consequence of the special form of the 
skein relation (2).   
However, the most difficult step is to prove that the canonical   
combinatorial extension from $K_{n}(\al,\be)$ to $K_{n+1}(\al,\be)$ 
is indeed a well-defined linear functional,   
which moreover satisfies the condition of trace commutativity.

The method of proof is greatly inspired by \cite{Berg}.    
One defines a graph whose vertices are linear combinations on    
elements of the free group in $n-1$ letters.  
The edges  correspond to pairs of elements which 
differ by exactly one relation, from the set of relations  
which present the algebras $K_{n}(\al,\be)$.     

Some of these edges will be given an orientation.
The first problem is whether each connected component 
of this graph has a minimal element for this orientation.  
We have to understand further whether different descending paths 
issued from the same vertex will eventually abut on 
the same element. 
Notice that whenever there is an unique minimal  
element in each component one is able to derive a basis    
for the module $K_{n}(\al,\be)$.  

In order to achieve the existence of minimal elements 
in each component one has to add a number of extra  
edges to our former graph. These new edges correspond 
to other relations  satisfied in $K_{n}(\al,\be)$.

Let us consider the lexicographic order on the letters   
generating the free group on $n-1$ letters. 
We want to use the relations in the algebra $K_n(\al,\be)$   
as transformations which  replace a word 
by a linear combination of  smaller ones. 
Using recursively this procedure the initial word 
is simplified until it reaches a normal form, where no more 
simplifications are possible.

The simplification procedure is encoded 
in the oriented paths of the graph: each relation used as above 
is an oriented edge of our graph.  
Specifically, these are  given by the 
following monomial substitutions:

 \begin{equation} 
{\rm (C0)(j)}: \;  a  b_j^3  b\rightarrow  \alpha \,  a  \bj^2   b  + \be\,  a  \bj  b +
a b
\end{equation}    
\begin{equation} 
{\rm (C1)(j)}:\; a b_{j+1}b_jb_{j+1} b \rightarrow  a b_jb_{j+1}b_j b 
 \end{equation}    
  \begin{equation} 
{\rm (C2)(j)}:\; a b_{j+1}b_j^2b_{j+1} b \rightarrow a S_j b
  \end{equation}    
   \begin{equation}
{\rm (C12)(j)}:\; a b_{j+1}b_j^2b_{j+1}^2 b \rightarrow a C_j b  
\end{equation}     
   \begin{equation}
{\rm (C21)(j)}:\; a b_{j+1}^2b_j^2b_{j+1} b \rightarrow a D_j b    
\end{equation} 
where  $E_{j+1}=\alpha  b_{j+1}^2   + \be  b_{j+1} + 1$, 
$S_j=  b_{j+1}b_j^2b_{j+1}-R_{0}(j)$ ,  $C_j=b_j^2b_{j+1}^2b_j+\alpha(b_{j+1}b_j^2b_{j+1}-b_jb_{j+1}^2b_j)+\beta(b_{j+1}b_j^2-b_{j+1}^2b_j)$ and    
$D_j=b_jb_{j+1}^2b_j^2+\alpha(b_{j+1}^2b_j^2b_{j+1}^2-b_jb_{j+1}^2b_j)+  
\beta(b_j^2b_{j+1}-b_jb_{j+1}^2)$, $j\in\{0, \dots, n-2\}$. 
Here $R_0(j)$ is the result of translating the indices of all
letters in $R_0$ by 
$j-1$ units.   
 
Several edges of our graph will remain  unoriented. 
The reason is that the respective relations 
are not compatible with the lexicographic order. 
They  correspond to the following monomial substitutions:    
   
\begin{equation}   
{(\rm P_{ij})}: \;  a b_ib_j b \rightarrow a b_jb_i b, \; \mbox{
  whenever } \;|i-j|>1 
\end{equation}  
The transformations (3)-(8) will be called reduction or 
simplification transformations.    
   
  Remark that we introduced some extra relations, 
namely (5) and (6), which  are not among the relations 
of the given presentation of $K_n(\al,\be)$, 
but which are nevertheless  satisfied in  
$K_n(\al,\be)$. This new relations  make the reduction process    
ambiguous. The reason for introducing them   
is to insure the existence of descending paths towards some    
minimal elements even in the case when the graph might  
contain closed oriented loops.

The next step consists of checking the existence and 
uniqueness of minimal  elements    
in this semi-oriented graph by means of 
so-called Pentagon Lemma (see section 2). 
One notices that one cannot always find a unique minimal element
by using directed paths issued from a fixed vertex.     
Furthermore we shall enlarge our graph to a tower of graphs    
modeling not one particular algebra $K_n(\al,\be)$ for fixed $n$,
but the set of  linear functionals defined  on the whole tower   
$\cup_{n=2}^{\infty}K_n(\al,\be)$  and satisfying     
certain compatibility conditions, which relate the values
taken on $K_n(\al,\be)$ to those on $K_{n+1}(\al,\be)$. 
The main feature of the tower is that now one can 
simplify further the minimal elements by recurrence 
on the level $n$, until one abuts on 
$K_0(\al,\be)$.  Here the Colored Pentagon Lemma 
(see section 3) can be applied and   
the uniqueness of the minimal elements in the tower of graphs 
is reduced to finitely many algebraic conditions.    
We will find  actually that the main 
obstructions lie in $K_4(\al,\be)$,   
as it might be inferred  from the study of 
quadratic Hecke algebras. From a different perspective, 
we actually proved that a certain linear functional on the tower 
$\cup_{n=2}^{\infty}K_n(\al,\be)$ is well-defined.

Eventually one has to verify whether the linear 
functional obtained above  satisfies the commutativity    
conditions for being a Markov trace. One proves that 
there is only one  obstruction to the commutativity, 
which lies  also in $K_{4}(\al,\be)$. 
   
Summarizing,  there are two types of  obstructions to the   
existence of a Markov traces:   
\begin{itemize}   
\item CPC obstructions, coming from the Colored Pentagon 
Condition,  and  
\item commutativity  obstructions.    
\end{itemize}    
These algebraic obstructions are polynomials with integer 
coefficients in the variables $\alpha$ and $\be$,  
and have been computed by using a computer code and 
formal calculus. The output of these computations 
is a set of explicit polynomials, which belong to the 
principal ideal generated by  $H_{(\alpha,\beta)}$. 
Furthermore,  the functional defined above is 
indeed a Markov trace, 
when restricting its values to the quotient by this 
principal ideal.

\section{Markov traces on $K_n(\al,\be)$}

\subsection{The cubic Hecke algebra $H(Q,3)$ revisited}   
   
The generalized Hecke algebras $H(P,3)$ could be considered 
for polynomials $P$ of higher degree by using the same 
definition as in the cubic case. One notices however
that $\dim_{\C}H(P,3)=\infty$ as soon as the degree 
of $P$ is at least $6$.    

\begin{remark}
The structure of the algebras $H(P,n)$ is well-known in 
the classical case (see \cite{Bou}) when $P$ is quadratic.    
They are finite dimensional semi-simple algebras of 
dimension $n!$, isomorphic (for generic $P$)
to the group algebra of the permutation group on $n$ elements. 
There is no general theory for higher degree polynomials $P$, 
due to their considerable complexity. 
\end{remark}
In the particular case of cubic $Q$ and $n=3$ 
it was shown in \cite{Fun} the following:    
   
\begin{proposition}   
If $Q$ is a cubic polynomial with $Q(0)\neq 0$ then      
 $\dim_{\bf \C}H(Q,3)=24$. 
A convenient base of the vector space $H(Q,3)$ is:   
  
$e_1=1, \, e_2=b_1, \, e_3=b_1^2, \, e_4=b_2, \, e_5=b_2^2, \, e_6=b_1b_2, \, e_7=b_2b_1, \, e_8=b_1^2b_2,$    
$ e_9=b_2b_1^2, \, e_{10}=b_1b_2^2, \, e_{11}=b_2^2b_1, \, e_{12}=b_1^2b_2^2, \,  e_{13}=b_2^2b_1^2, \, e_{14}=b_1b_2b_1,$   
$e_{15}=b_1^2b_2b_1,  \, e_{16}=b_1b_2b_1^2, \, e_{17}=b_1b_2^2b_1^2, \, e_{18}=b_1^2b_2b_1^2, \, $   
$e_{19}=b_1^2b_2^2b_1, \, e_{20}=b_1b_2^2b_1,   \,    
e_{21}=b_1^2b_2^2b_1^2,  \,    
e_{22}=b_2b_1^2b_2, \, $   
$e_{23}=b_2b_1^2b_2b_1=b_1b_2b_1^2b_2, \, e_{24}=b_2b_1^2b_2b_1^2=b_1b_2b_1^2b_2b_1=b_1^2b_2b_1^2b_2$ 
\end{proposition}

\begin{proposition}  
$H(Q,3)$ is a semi-simple algebra which decomposes  generically as   
$   
\C^3 \oplus M_2^{\oplus3} \oplus M_3   
$,  where  $M_n$ is the algebra of $n \times n$ matrices. The morphism   into $H(Q,3)\to \C^3$ is obtained via the 
abelianization map. Each one of the three 
projections $H(Q,3)\to M_2$ factors through the projection 
$H(Q,3)\to H(P_i,3)=\C^2 \oplus M_2$ onto 
the quadratic Hecke algebra $H(P_i)$ defined by one divisor
$P_i$ of $Q$.  
\end{proposition}
\begin{proof}   
This follows by  a direct computation, making use of the 
following identities (\cite{Fun}):

\[b_{j+1}b_j^2b_{j+1}b_j=b_jb_{j+1}b_j^2b_{j+1}\]   
\[ b_{j+1}^2b_j^2b_{j+1}=b_jb_{j+1}^2b_j^2+\alpha(b_{j+1}b_j^2b_{j+1}-b_jb_{j+1}^2b_j)+ \beta(b_j^2b_{j+1}-b_jb_{j+1}^2)\]   
\[ b_{j+1}b_j^2b_{j+1}^2=b_j^2b_{j+1}^2b_j+
\alpha(b_{j+1}b_j^2b_{j+1}-b_jb_{j+1}^2b_j)+ 
\beta(b_{j+1}b_j^2-b_{j+1}^2b_j)\]   
\end{proof}

\subsection{The  algebras $K_n(\al,\be)$}

The tower $\cup_{k=1}^{\infty}P(k)$ of  quotients 
of $\cup_{k=1}^{\infty}H(Q,k)$ is homogeneous if any    
identity $F(b_i,b_{i+1},\dots,b_j)=0$    
which holds in $P(j+1)$, remains valid under the translation 
of indices i.e. 
$F(b_{i+k},b_{i+k+1},\dots,b_{j+k})=0$, for  all 
$k\in \Z$ such that  $k\geq 1-i$.
If one seeks for Markov traces on towers of quotients of 
$\C[B_n]$ it is convenient to restrict ourselves to the study
of homogeneous quotients.

We define $K_n(\al,\be)$ as the homogeneous quotient
$\hqn / I_n$, where $I_n$ is the two-sided   
ideal  generated by:   
       
 $   
\bj  \bjmu^2  \bj +   
(\be^2 - \al) \bjmu^2 \bj^2 \bjmu^2  +   
 (\al^2 - \al \be^2 - \be) \bjmu \bj^2  \bjmu^2 +   
$   
$   
(\al^2 - \al \be^2 - \be) \bjmu^2 \bj^2    
$   
$   
 \bjmu   
+ (\al^2 - \al \be^2) \bjmu^2 \bj \bjmu^2   
+   
$   
$   
(1+ 2\al \be  + \al^2 \be^2 - \al^3)   \bjmu   \bj^2  \bjmu +     
(1+ \al \be  + \al^2 \be^2 -   
$   
$   
 \al^3) \bjmu   \bj   \bjmu^2  +   
(1+ \al \be  + \al^2 \be^2 - \al^3) \bjmu^2    \bj   \bjmu   
+(1+ 2\al \be  - \be^3)\bj^2   \bjmu^2  +     
$   
$   
 (1+ 2\al \be  - \be^3)\bjmu^2   \bj^2  +     
(\al \be^3 - 2 \al - 2 \al^2 \be) \bj   \bjmu^2   
+   
$   
$   
(\al \be^3 - 2 \al - 2 \al^2 \be) \bjmu^2   \bj +   
$   
$   
 (\al \be^3 - 2 \al - 2 \al^2 \be + \be^2) \bj^2   \bjmu    
+(\al \be^3 - 2 \al -   
 2 \al^2 \be + \be^2) \bjmu   \bj^2 + (\al^4 -\al^3 \be^2 -   
 $   
 $   
 2\al^2\be - 3 \al) \bjmu \bj \bjmu +   
(2 \al^3 \be + 3\al^2 - \al^2 \be^3 - \al \be^2) \bj \bjmu +   
$   
$   
(2 \al^3 \be + 3\al^2 -   
$   
$   
\al^2 \be^3 -   
 \al \be^2) \bjmu \bj +   
(\be^4-2\be-   
3\al\be^2+\al^2) \bjmu^2 +   
$   
$   
(\be^4-2\be-3\al\be^2+\al^2) \bj^2 +   
$   
$   
 (1+ 4\al \be +   
 3\al^2 \be^2 - \al^3 - \al \be^4- \be^3)   
\bjmu +   
(1+ 3\al \be + 3\al^2 \be^2 - \al^3 - \al \be^4)\bj   
$   
$   
 + 3\be^2 - \be^5 - 2 \al - 3 \al^2 \be + 4 \al \be^3    
$   

\noindent for $j\in \{1,\dots,n-1\}$.

\begin{proposition}  
Under the identification $H(Q,3)\cong     
\C^3 \oplus M_2^{\oplus3} \oplus M_3   
$,  the quotient
$K_{3}(\al,\be)$ corresponds to $M_2^{\oplus3} \oplus M_3$.   
\end{proposition} 
\begin{proof} 
In fact, it  suffices to show that the ideal    
$I_3$ is a vector space of dimension $3$.   
Let $I$ be the span of $R_0, R_1, R_2$, where:   
      
  $   
R_0 = \bm   \bc^2 \bm +   
 (\be^2 - \al) \, \bc^2  \bm^2  \bc^2  +  $   $   
 (\al^2 - \al \be^2 - \be)  \bc \bm^2  \bc^2 +   
(\al^2 - \al \be^2 - \be)  \bc^2 \bm^2   \bc^2    $  $    
+ (\al^2 - \al \be^2)  \bc^2  \bm  \bc^2 +   
$ $   
(1+ 2\al \be  + \al^2 \be^2 - \al^3)    \bc    \bm^2    \bc +     
(1+ \al \be  + \al^2 \be^2 - \al^3)  \bc    \bm    \bc^2  +  $   $   
(1+ \al \be  + \al^2 \be^2 - \al^3)  \bc^2     \bm    \bc   
$ $   
+(1+ 2\al \be  - \be^3) \bm^2    \bc^2  +     
 (1+ 2\al \be  - \be^3) \bc^2    \bm^2  +    $ $   
 (\al \be^3 - 2 \al - 2 \al^2 \be)  \bm    \bc^2   
$ $   
 +(\al \be^3 - 2 \al - 2 \al^2 \be)  \bc^2    \bm +   
 (\al \be^3 - 2 \al - 2 \al^2 \be + \be^2)  \bm^2    \bc    $  $   
+(\al \be^3 - 2 \al - 2 \al^2 \be + \be^2)  \bc    \bm^2 +   
 $ $   
(\al^4 -\al^3 \be^2 -2\al^2\be - 3 \al)  \bc  \bm  \bc +   
(2 \al^3 \be + 3\al^2 - \al^2 \be^3 - \al \be^2)  \bm  \bc +   $ $   
(2 \al^3 \be + 3\al^2 - \al^2 \be^3 - \al \be^2)  \bc  \bm +   
$ $    
(\be^4-2\be-3\al\be^2+\al^2)  \bc^2 +   
(\be^4-2\be-3\al\be^2+\al^2)  \bm^2    $ $   
 + (1+ 4\al \be + 3\al^2 \be^2 - \al^3 - \al \be^4- \be^3)  \bc $ $   
+ (1+ 3\al \be + 3\al^2 \be^2 - \al^3 - \al \be^4)  \bm   
 + 3\be^2 - \be^5 - 2 \al - 3 \al^2 \be + 4 \al \be^3    
$

\vspace{10pt}

  $   
R_1 = \bc  R_0 = \bc  \bm   \bc^2   \bm    
-\be  \bc^2  \bm^2  \bc^2  +   
 (1 + \al \be)  \bc  \bm^2   \bc^2 +   
 (1 + \al \be)  \bc^2  \bm^2   \bc^2  $ $   
+(1 + \al \be)  \bc^2  \bm  \bc^2   
$ $   
(-\al^2 \be- 2\al)  \bc   \bm^2    \bc +     
(-\al^2 \be- 2\al)  \bc   \bm    \bc^2 +  $ $   
(-\al^2 \be- 2\al)  \bc^2   \bm    \bc +   
$ $   
(\be^2 - \al) \bm^2    \bc^2  +     
(\be^2 - \al) \bc^2    \bm^2  +    $ $   
(\al^2 - \al\be^2)  \bm    \bc^2 +   
$ $   
(\al^2 - \al\be^2)  \bc^2    \bm +   
(\al^2 - \al\be^2 - \be)  \bm^2    \bc +  $ $   
(\al^2 - \al\be^2 - \be)  \bc    \bm^2 +   
$ $   
(\al^3 \be + \be  + 3\al^2)  \bc  \bm  \bc +   
(1+ \al \be  + \al^2 \be^2 - \al^3)  \bm  \bc +  $ $   
(1+ \al \be  + \al^2 \be^2 - \al^3)  \bc  \bm +   
$ $   
(1+ 2\al \be  - \be^3)  \bc^2 +   
(1+ 2\al \be  - \be^3)  \bm^2 +   
(\al \be^3 - 2 \al - 2 \al^2 \be + \be^2)  \bc   
$ $   
+ (\al \be^3 - 2 \al - 2 \al^2 \be )  \bm $  $   
+\be^4-2\be-3\al\be^2+\al^2 
$

\vspace{10pt}

  $   
R_2 = \bc  R_1 =  \bc^2  \bm   \bc^2   \bm  +$ $    
 \bc^2  \bm^2  \bc^2     
-\al  \bc  \bm^2   \bc^2    
-\al  \bc^2  \bm^2   \bc^2  $ $   
-\al  \bc^2  \bm  \bc^2+   
$ $   
\al^2  \bc   \bm^2    \bc +     
(\al^2 + \be)  \bc   \bm    \bc^2 +   
(\al^2 + \be)  \bc^2   \bm    \bc + $ $   
 $ $   
(- \be)   \bm^2    \bc^2  +     
(- \be)   \bc^2    \bm^2  +   
(1 + \al \be)  \bm    \bc^2 +   
 $ $   
(1 + \al \be)  \bc^2    \bm +   
(1 + \al \be)  \bm^2    \bc + $ $   
(1 + \al \be)  \bc    \bm^2 +   
 $ $   
(-\al^3 \be -\al \be  + 1)  \bc  \bm  \bc +   
(-\al^2 \be- 2\al)  \bm  \bc + $ $    
(-\al^2 \be- 2\al)  \bc  \bm +   
 $ $   
(\be^2 - \al)  \bc^2 +   
(\be^2 - \al)  \bm^2 +   
(-\al \be^2+ \al^2 - \be)  \bc +   
 $ $   
(-\al \be^2+ \al^2)  \bm +   
 1+ 2\al \be  - \be^3    
 $

\begin{lemma}   
There is an isomorphism of vector spaces $I\cong I_3$.   
\end{lemma}   
\begin{proof}   
Remark first that the following identities hold true in $H(Q,3)$:    
$$   
b_1R_0=R_0b_1=R_1, \; b_1R_1=R_1b_1=R_2, \; b_1R_2=R_2b_1=R_0 + \be R_1+ \al R_2 
$$    
Then, by direct computation, we obtain that:    
$$b_2R_0=R_0b_2=R_1, \;  b_2R_1=R_1b_2=R_2, \; b_2R_2=R_2b_2=R_0 + \be R_1+ \al R_2$$    
From these relations we derive that 
$xR_0y\in I$ for all $x,y \in H(Q,3)$, and hence 
$I_3 \subset I$.  The other  inclusion is  immediate.   
\end{proof}   

\noindent The proposition is then a consequence of the previous 
lemma.\end{proof}

\subsection{Uniqueness of the Markov trace on the tower 
$\cup_{n=1}^{\infty}K_n(\al,\be)$}

From now on we will work  with the group ring $\Z\left[\al, \be\right]\left[B_{n}\right]$ instead of    
$\C \left[B_{n}\right]$.

\begin{definition}    
Let $z,\bar{z}\in \Z(\al,\be)$ be rational functions 
in the variables $\al$ and $\be$, and $R$ a 
$\Z\left[\al, \be, z,\bar{z} \right]$-module.
The linear functional 
$\mct: \cup_{n=1}^{\infty} K_n(\al,\be)\to R$ is 
said to be an admissible functional (with parameters $z$ and 
$\bar{z}$) on $\cup_{n=1}^{\infty}K_n(\al, \, \be)$ if the 
following conditions are fulfilled:   
\[ \mct(xb_ny)=z\mct(xy) \; \; \mbox{ for all} \; \;  x,y\in 
K_n(\al,\be)\]    
\[ \mct(xb_n^{-1}y)=\bar{z}\mct(xy) \; \; \mbox{ for  all} \; \; x,y\in  K_n(\al,\be)\]

\noindent An admissible functional $\mct$   
is a Markov trace if it satisfies the following trace condition:    
\[     
\mct(a b)=\mct(b a) \;\; \mbox{for  any}  \;\; a, b \in  
K_n(\al,\be)
\]   
\end{definition}

\begin{remark}
The tower of quadratic Hecke algebras admits an unique 
Markov trace (\cite{Jones}). Similarly, the tower of 
Birman-Wenzl algebras (\cite{BW}) admits an unique Markov trace. 
\end{remark}  
 \begin{definition}  
The  admissible functional  $\mct$ is {multiplicative}     
if $\mct(xb_n^k)=\mct(x)\mct(b_n^k)$  holds for all $x\in H(Q,n)$
and $k\in \Z$.   
\end{definition}   
\begin{remark} The Markov trace on the quadratic Hecke algebras
is multiplicative, and hence     
$\mct(xy)=\mct(x)\mct(y)$ for any     
$x\in H(Q,n)$ and 
$y\in \langle 1,b_n,b_{n+1},\dots,b_{n+k}\rangle$.   
However, one cannot expect that this property holds true 
for Markov traces on arbitrary higher degree Hecke   
algebras. 
\end{remark}

\begin{proposition}  
The admissible functionals on the tower of cubic Hecke 
algebras are multiplicative. In particular: 
\[\mct(a b_n^2 b)= t \mct(ab) \; \mbox{ for all } a,b \in H(Q,n) \]  
where $t = \al z + \be +\bar{z}$.
\end{proposition}
\begin{proof} One uses the identity 
$b_n^2=\alpha b_n+\beta+ b_n^{-1}$ for proving 
the multiplicativity for $k=2$, and then 
continue by recurrence for all $k$. 
\end{proof}

\noindent One  can state now the unique extension property of Markov traces.

\begin{proposition}    
For fixed $(z,t)$ there exists at most one    
Markov trace on $K_n(\al,\be)$ with parameters $(z,t)$.    
\end{proposition}    
\begin{proof} Define recursively the modules $L_n$ as follows:

\[L_2=H(Q,2)\]    
\[L_3=\C\langle b_1^ib_2^jb_1^k| \mbox{ where } 
i, j, k \in \{0,1,2\}\rangle\]    
\[L_{n+1}=\C \langle ab_n^{\varepsilon}b \, | \mbox{ where }
\, a, b \mbox{ are elements of the basis of} \; L_n , \,
\mbox{ and } \varepsilon\in \{1,2\}\rangle\; \oplus \; L_n\]

\begin{lemma}\label{sim} 
The natural projection  
$\pi:L_n\to K_n(\al,\be)$ is surjective.   
\end{lemma}   
\begin{proof} For    
$n=2$ it is clear. For $n=3$ we know that $ \bm   \bc^2   \bm, \; \bc  \bm  \bc^2  \bm , \; \bc^2     
\bm      
\bc^2  \bm \in $ $\pi(L_3)$, from the exact form of the relations
$R_0, R_1, R_2$, generating the ideal $I_3$. 
We shall use a recurrence on $n$ and assume that the 
claim holds true for $n$. 
   
Consider now $w\in K_{n+1}(\al,\be)$ represented by a word in the  
$b_i$'s having only positive exponents. We assume that the    
degree of the word in the variable $b_n$ is minimal among all  
linear combinations of  words (with positive exponents) 
representing $w$.

\begin{enumerate}   
\item If this degree is less or equal to 1 then there is nothing to prove.    
  
\item If the degree is 2 then either $w=ub_{n}^2v$, $u,v\in  K_n(\al,\be)$ so using the induction hypothesis we are done, or else    
$w=ub_nzb_nv$, where $u,z,v\in  K_n(\al,\be)$. Therefore $z=xb_{n-1}^{\varepsilon}y$ where $x,y\in    
K_{n-1}(\al,\be)$ by the induction hypothesis and    
$\varepsilon\in\{0,1,2\}$. 
\begin{enumerate}
\item If $\varepsilon = 0$ then $w$ can be reduced to $uzb_n^2v$. 
\item If $\varepsilon=1$ then    
$w=ub_nxb_{n-1}yb_nv=uxb_{n-1}b_nb_{n-1}yv$ 
hence the degree of $w$  can be lowered by one, 
which  contradicts our minimality assumption.   
\item If $\varepsilon=2$ then    
$w=uxb_nb_{n-1}^2b_nyv$.  One derives that:      
\[b_nb_{n-1}^2b_n\in \C \langle b_{n-1}^ib_n^jb_{n-1}^k,\; i,j,k\in \{0,1,2\}\rangle\]    
hence we reduced the problem to the case when   
$w$ is a word of type $u'b_n^2v'$.    
\end{enumerate}  
    
\item If the degree of $w$ is at least 3 we will 
contradict the minimality assumption.    
In fact, in this situation $w$ will contain either a sub-word   
$w'=b_n^aub_n^b$, with $u\in K_n(\al,\be)$ 
and $a+b\geq 3$, or else    
a sub-word $w''=b_nub_nvb_n$, with $u,v\in K_n(\al,\be)$.    
   
\begin{enumerate}    
\item In the first case using the induction we can write    
 $u=xb_{n-1}^{\varepsilon}y$, with $x,y\in K_{n-1}(\al,\be)$.    
   
\begin{enumerate}    
\item Furthermore, if $\varepsilon=0$ then $w'=b_n^{a+b}xy=\alpha
b_n^{a+b-1}xy+\beta b_n^{a+b-2}xy+b_n^{a+b-3}xy$, and   
hence the degree of $w$ can be lowered by one.

\item If $\varepsilon=1$ then    
$w'=b_{n}^{a-1}xb_nb_{n-1}b_nyb_n^{b-1}=b_n^{a-1}xb_{n-1}b_nb_{n-1}yb_n^{b-1}$,   
and again its degree can be reduced by one unit.

\item If $\varepsilon=2$ then either $a$ or $b$ is  equal 2. 
Assume that $a=2$. We can  therefore write:     
\begin{eqnarray*}
w' & = &
xb_n^2b_{n-1}^2b_nyb_{n}^{b-1}=xb_{n-1}b_n^2b_{n-1}^2yb_n^{b-1}+ 
\alpha(b_nb_{n-1}^2b_n-b_{n-1}b_n^2b_{n-1})yb_n^{b-1}+ \\
& &  \beta(b_{n-1}^2b_{n}-b_{n-1}b_n^2)yb_n^{b-1}
\end{eqnarray*}    
contradicting again the minimality of the degree of $w$.    
\end{enumerate}

\item In the second case we can write also $u=xb_{n-1}^{\varepsilon}y$, $v=rb_{n-1}^{\delta}s$ with $x,y,r,s\in K_{n-1}(\al,\be)$.

\begin{enumerate}   
\item  
If $\varepsilon$ or $\delta$ equals 1 then, after some obvious  
 commutations the word $w"$ contains the sub-word     
$b_nb_{n-1}b_n$ which can be replaced by 
$b_{n-1}b_nb_{n-1}$ and hence  diminishing its degree.

 \item  If $\varepsilon=\delta=2$ then    
$w"=xb_nb_{n-1}^2b_nyrb_{n-1}^2b_ns$.   
We use the homogeneity to replace $b_nb_{n-1}^2b_n$ by a sum of elements of type $b_{n-1}^ib_n^jb_{n-1}^k$.    
 Each term of the expression of $w"$ which comes from a factor 
which has the exponent $j < 2$, has diminished its degree. 
 The remaining terms are     
 $xb_{n-1}^ib_n^2b_{n-1}^kyrb_{n-1}^2b_ns$,    
 so they contains a sub-word $b_n^2ub_n$ whose degree we already know  
 that it can be reduced as above. This proves our claim. 
\end{enumerate}
\end{enumerate}
\end{enumerate}
\end{proof}

\noindent Eventually recall that the Markov traces $\mct$ on 
$\cup_{n=1}^{\infty}H(Q,n)$ are multiplicative,
and hence they satisfy:    
 $\mct(xb_n^{\varepsilon}y)=\mct(b_n^{\varepsilon})\mct(yx)$.    
Therefore there is a unique extension of $\mct$ from 
$K_{n}(\al,\be)$ to  $K_{n+1}(\al,\be)$. 
This ends the proof of our  proposition. 
\end{proof}

\begin{proposition}\label{indu}
The admissible functionals on the tower of 
algebras $\cup_{n=1}^{\infty}K_n(\al,\be)$ satisfy the identities: 
\[\mct(xuv)=\mct(u)\mct(xv)\; \mbox{ for } x,v\in H(Q,m)
\;\mbox{  and } u\in \langle 1,b_m,b_{m+1},\dots,b_{m+k}\rangle\] 
\end{proposition}
\begin{proof}
For $k=0$ this is equivalent to the multiplicativity 
of the admissible functional. 
We will use a recurrence on $k$ and assume that 
the claim holds true for  $k$. 
By lemma \ref{sim} one can reduce the element $u$ 
in $K_{m+k+1}(\alpha,\beta)$ to a 
(non-necessarily unique)  normal form 
$u=u_1b_{m+k}^{\varepsilon}u_2$, where 
$u_j\in \langle 1,b_m,b_{m+1},\dots,b_{m+k}\rangle$, $j\in\{1,2\}$  
and $\varepsilon\in\{0,1,2\}$.  The multiplicativity of the 
admissible functionals implies that: 
\[\mct(xuv)=\mct(b_{m+k}^{\varepsilon})\mct(xu_1u_2v)\] 
By the recurrence hypothesis one knows that: 
\[ \mct(xu_1u_2v)= \mct(u_1u_2)\mct(xv) \]
and since:
\[ \mct(u)=\mct(b_{m+k}^{\varepsilon})\mct(u_1u_2) \]
one derives our claim. \end{proof}

\section{CPC Obstructions}

\subsection{The pentagonal condition}

The following  is an immediate 
consequence of lemma \ref{sim}:    
   
\begin{lemma} \label{lem:man}    
There is a surjection of $(K_n(\al,\be),K_n(\al,\be))$-bimodules:    
\[ K_n(\al,\be)\oplus K_n(\al,\be)\otimes_{K_{n-1}(\al,\be)}K_n(\al,\be) \oplus K_n(\al,\be)\otimes_{K_{n-1}(\al,\be)}K_n(\al,\be)\longrightarrow K_{n+1}(\al,\be)\]    
 given by:    
\[x\oplus y\otimes z \oplus u\otimes v\rightarrow x+yb_nz+ub_n^2v
\]    
\end{lemma}

\begin{remark}   
In particular, the admissible functionals on the tower 
$\cup_{n=1}^{\infty}K_n(\al,\be)$  are unique up to    
the choice of $\mct(1)\in R$.    
\end{remark}
Now, we want to use the transformations (3)-(7) to simplify 
the positive words from $K_n(\al,\be)$, so that 
the degree of $b_{n-1}$  becomes as small as possible. 
According to   
the previous lemma every word in $K_n(\al,\be)$ can be written as 
a linear combination of  words of the form 
$x_ib_{n-1}^{\varepsilon_i}y_i$, with $\varepsilon_i\in\{0,1,2\}$
and $x_i,y_i\in K_{n-1}(\al,\be)$. 
Unfortunately, one needs to use in both directions 
the transformations ${\rm P}_{ij}$ from (8): \;     
$b_ib_j\leftrightarrow b_jb_i, \;\mbox{ for } \; | i-j| >1$.  
  
\begin{remark}
The linear combination we obtained above is 
a kind of  {\em normal form} for the    
word with which we started.  
It could happen that this normal form is not 
unique since we may perform again permutations of type (8) 
among some of its letters. 
However, if any two such normal forms were    
equivalent under  the transformations (8), then we would 
obtain an  almost canonical description  
of the basis of $K_n(\al,\be)$. 
This assumption is equivalent to saying that the  
surjection from lemma \ref{lem:man} is an isomorphism.   
Unfortunately, this is not the case. 
However, one can describe the obstructions   
to the uniqueness for this almost canonical form, as follows.   
\end{remark}  
We return now to the module of the   
admissible functionals on the whole tower of algebras 
$\cup_{n=1}^{\infty}K_n(\al,\be)$. The conditions satisfied by admissible
functionals enable us to add a new type of 
simplifications, 
by means of the following formulas:   
   
\begin{equation}
ab_{n-1}b\rightarrow z \:ab, \; \mbox{ and respectively }\; 
ab_{n-1}^2b\rightarrow t \: ab,
\; \mbox{ where } a,b\in K_{n-1}(\al,\be)
\end{equation}   
This way we can reduce a word from $K_n(\al,\be)$ 
to a linear combination  of words from 
$K_{n-1}(\al,\be)$. Assume that we are using repeatedly
the transformations (9). Then we will eventually reduce 
the initial word to a linear combinations of words 
in $K_0(\al,\be)$, thus to an element of $R$.  
Remark that this element is actually 
the value that the admissible functional takes on the 
initial word.   
Our main task is to understand whether the
final reduction is independent on the way we chose 
to make the simplifications. When this happens to be true
then we obtain that the  functional which associates to 
each element of $K_n(\al,\be)$ its final reduction
is a well-defined  admissible functional. 
However, we will encounter  below some obstructions to the 
uniqueness, which fortunately we can treat explicitly.

One formalizes this procedure at follows.    
Let $\Gamma$ be a semi-oriented graph. This means that 
some of its  edges are oriented while the remaining 
ones are left unoriented.  We write  $v\rightarrow w$
if there is an oriented edge from $v$ to $w$.   
A path $v_1v_2\cdots v_n$ in $\Gamma$ is called a 
{\em semi-oriented path} if, for each 
$j$, one has  either  $v_j\rightarrow v_{j+1}$ or else   
$v_jv_{j+1}$ is an  unoriented edge of $\Gamma$.    
If all edges of the path  are unoriented then 
we say that its endpoints are  {\em (weakly) equivalent}.    
   
\begin{definition} 
The sequence of vertices    
$[v_0,v_1,\dots ,v_{n+1}]$ is an {\em open pentagon configuration} in 
$\Gamma$ (abbreviated o.p.c.) if $v_1\rightarrow v_0$, $v_1v_2\cdots v_{n-1}$    
is an unoriented path and  $v_{n}\rightarrow v_{n+1}$.   
\end{definition}

\begin{definition}   
The semi-oriented graph 
$\Gamma$ verifies the pentagon condition (abbreviated PC) 
if for any  open pentagon configuration  
$[v_0,v_1,\dots,v_{n+1}]$  there exist    
semi-oriented paths $v_0x_1x_2\cdots x_me$ and $v_{n+1}y_1y_2\cdots y_pe$ having the    
same endpoint. 
\end{definition}

Given a graph like above one has a binary relation 
induced as follows:  we set $x\leq y$ if there exists an semi-oriented path from $y$ to $x$ in $\Gamma$. Of course $\leq$ is not    
always a partial order relation. 
A necessary and sufficient condition  for $\leq $ to be 
a partial order is  that $\Gamma$ contains no closed 
semi-oriented  closed loops.    
One says  that $x$ is {\em minimal} if $y\leq x$ 
implies that $y$ is weakly equivalent to $x$.

\begin{lemma} \label{lem:1}   
Suppose that the (PC) holds. If a  connected 
component $C$ of the graph $\Gamma$ has a    
minimal element then this is unique up to weak 
equivalence.   
\end{lemma}    
\begin{proof} Consider two minimal  elements $x$ and $y$ which lie in $C$. Then there exists some path    
$xx_0x_1\cdots x_ny$ joining  them. Since $x$ is minimal the closest oriented    
edge  - if it exists - must be in-going; and the same    
is true for $y$. If this  path is not unoriented, then  the
minimality  implies that there are at least two oriented 
edges. Therefore one 
can find a sequence of  open pentagon configurations lining 
on the path which joins $x$ to $y$. 
We apply then  the (PC) iteratively, whenever we see one    
such o.p.c.,  or one o.p.c.  appears at the next stage, as 
in the figure below:

\begin{center}
\hspace{10pt}\psfig{figure=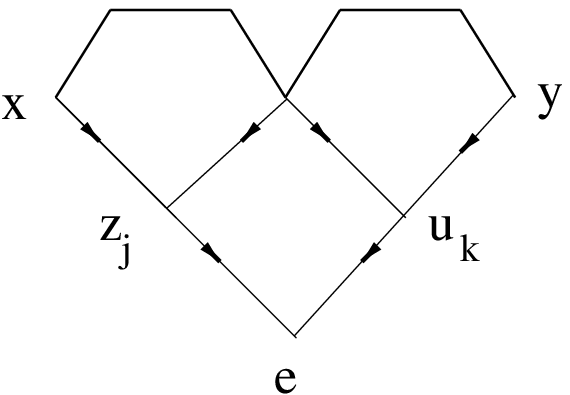,width=4cm}   
\end{center}

\noindent When this process stops, we find two semi-oriented 
paths $xz_1z_2\cdots z_pe$ and $yu_1u_2\cdots u_se$    
having the same endpoint $e$. So $e\leq x$ and $e\leq y$. 
From  minimality both these paths must be unoriented, and thus 
$x$ and $y$ are  weakly equivalent. \end{proof}

\begin{remark}  
The existence of minimal elements is not {\em a priori} granted, 
without additional conditions. If $\leq$ had been  
 a partial order with descending chain condition, 
then the existence of  
 minimal elements would be  standard. We will show that in the 
present case, of the graph modeling the admissible 
functionals on the tower $\cup_{n=1}^{\infty}K_n(\al,\be)$,  such 
minimal elements exist, though as $\leq$ is not a partial order. 
\end{remark}

\subsection{The colored  tower of graphs $\Gamma^*_n$}

Suppose now that  we have a family of semi-oriented graphs $\Gamma_n$
as follows.  Each graph $\Gamma_n$ has 
a distinguished subset of vertices $V_n^0$ 
whose elements are minimal elements in their connected    
respective components.  Assume also that each connected 
component of the $\Gamma_n$ admits at least 
one minimal element.  Further, we suppose that 
each vertex from $V_n^0$ has exactly one outgoing 
edge which joins it to a vertex of $\Gamma_{n-1}$. 
We color the edges connecting  graphs  
$\Gamma_n$ and $\Gamma_{n-1}$ in red. 
Set $\Gamma_n^*$ for the union of all $\Gamma_j$, with $j\leq n$
to which we add all red edges connecting graphs 
$\Gamma_k$ and $\Gamma_{k+1}$, for $k\leq n-1$.
We can have an intuitive view of $\Gamma_n^*$ by looking at the 
$\Gamma_n$ as graphs lying on different floors which are 
connected by vertical red edges pointing downwards. 
   
\begin{definition}   
The graph $\Gamma_n^*$ is {\em coherent} if any connected component   
of $\Gamma_n$ has an   unique minimal element within $\Gamma_n^*$, 
 up    to weak equivalence.
\end{definition} 
\begin{remark} A minimal element should belong to $\Gamma_0$. 
\end{remark} 
We  state now the {\em colored} version of the Pentagon Lemma for this type of graphs.

\begin{definition}   
We say that $\Gamma_n$ 
verifies the colored pentagon condition  (CPC)  
if, for any open pentagon configuration 
$[v_0,v_1,\dots,v_{m+1}]$  in $\Gamma_n$,     
 there exist bicolored semi-oriented paths (in $\Gamma_n^*$) 
from $v_0$ and $v_{m+1}$ having the same endpoint.    
 In addition, if $xy$ is an unoriented edge in $\Gamma_n$ with $x,y\in V_n^0$ then there exist    
 semi-oriented paths in $\Gamma_n^*$ starting with red edges and
 having the same endpoint, as in the figure below:  
    
\hspace{72pt}\psfig{figure=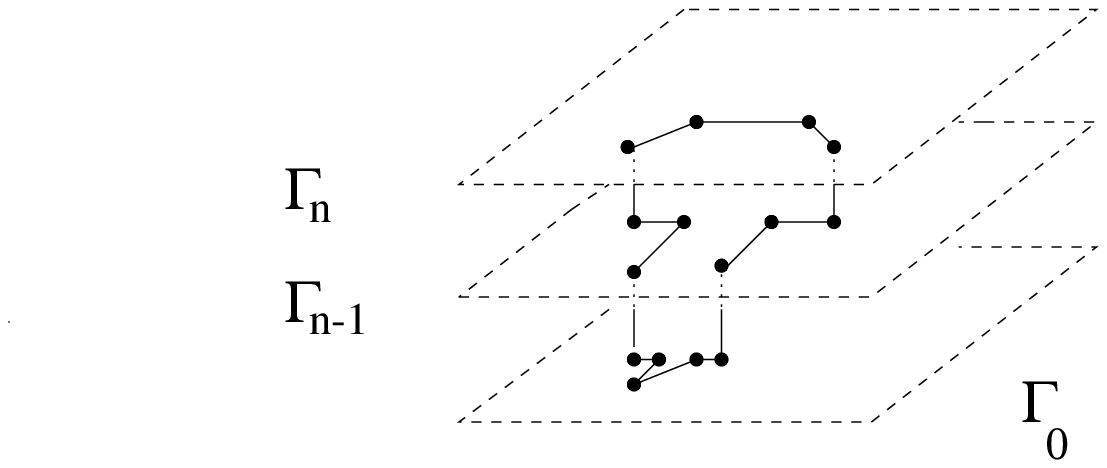,width=6.5cm}   
\end{definition}

\begin{lemma} \label{lem:2}   
Suppose that $\Gamma_{n-1}^*$ is coherent and the (CPC) condition is fulfilled. Then    
$\Gamma_n^*$ is coherent.    
\end{lemma}    
\begin{proof} The proof is similar to that of Pentagon Lemma. \end{proof}

Now, we are ready to define the sequence of semi-oriented 
graphs $\Gamma_n$, which models the admissible functionals 
on $\cup_{n=1}^{\infty}K_n(\al,\be)$.
\begin{definition}  
The vertices of $\Gamma_n$  are the elements of the ring algebra   
$\Z [ \al, \be, z,  \bar{z}][F_n]$, where $F_n$ is    
the free monoid $F_{n-1}$ generated by $n-1$ letters 
$\{b_1,b_2,\dots,b_{n-1}\}$.  
The vertices of  $\Gamma_0$ are the elements of    
$\Z [ \al, \be, z,  \bar{z}]$.  
Two vertices    
$v=\sum_{i}\alpha_i\; x_i$ and $w=\sum_i\beta_i\; y_i$, 
where $\alpha_i,\beta_i\in
\Z [\al ,\be, z,\bar{z}]$  and $x_i, \, y_i \in F_n$, 
are related    
by an oriented edge if exactly one monomial $x_i$  of $v$ is 
changed by means of a reduction transformation
among the rules (3)-(7).   
An unoriented edge between $v$ and $w$ corresponds to a 
simplification transformation (8) of one 
monomial $x_i$ from the previous expression of  $v$.   
\end{definition}  
\begin{remark}    
The use of (C12) and (C21) is somewhat ambiguous 
since we can always use (C2) for a    
sub-word of the given word. Their role is to break in some sense the
closed oriented loops in $\Gamma_n$, as we shall see below.
\end{remark}
Consider now the following sets of words in the $b_i$'s:         
\[W_0=\{1\}\]    
\[W_{n+1}=W_n\cup W_nb_{n+1}Z_n\cup W_nb_{n+1}^2Z_n \]    
where:    
\[Z_n=\{b_{n}^{i_0}b_{n-1}^{i_1}\cdots b_{n-p}^{i_p}| \mbox{ where 
the indices }
 i_1,i_2,\dots,i_p\in\{1,2\}, \mbox{ and } p\in\{0,1,\dots,n-1\}\}\]
Let $V_n^0$ be the set of vertices corresponding  to elements 
of the $\Z [ \al ,\be, z,\bar{z}] $-module    
generated by $W_n$. This completes the definition of the tower of
graphs $\Gamma_n$. We have the following result:    
\begin{proposition} \label{lem:3}   
Each connected component of $\Gamma_n$ has a minimal 
element in $V_n^0$, not necessarily unique.    
\end{proposition}    
\begin{proof} We use an induction on $n$. 
For $n=0$ the claim is obvious. 
Let now $w$ be a word in the $b_i$'s having 
only positive exponents. 
\begin{enumerate}
\item If its degree in $b_{n}$ is zero or one, then we apply the 
induction  hypothesis and we are done. 
\item If the degree in $b_n$ is 2 and $w$ contains the   
sub-word  $b_{n}^2$, then again we are able to apply 
the induction hypothesis.   
\item By using (C0) several times one  can also  
suppose that no exponents greater than 2 occur in $w$. 
\begin{enumerate}
\item If the degree of $b_n$ is 2 then $w= x b_{n}yb_{n}z$ with 
$x,y,z\in F_{n-1}$. The induction hypothesis applied 
to $y$ implies that $w\geq xb_{n}ab_{n-1}^{\varepsilon}bz$ 
with $a,b\in F_{n-1}$. Then several transforms of 
type ($P_{n j}$) and (C$\varepsilon$) will do the job.    
\item  Consider now that the degree in $b_n$ is at least 3.    
Then $w$ contains a sub-word which has either the form
$b_{n}^{\alpha}xb_{n}^{\beta}$ with $3\leq\alpha+\beta\leq 4$, 
or else one of the type $b_{n}xb_{n}yb_{n}$.    
The second case reduce to the first one as above.    
In the first case assume that $x\geq ab_{n-1}^{\varepsilon}b$ 
for some $a,b\in F_{n-2}$.  Then several applications 
of ($P_{n j}$) lead us to consider the 
sub-word $b_{n}^{\alpha}b_{n-1}^{\varepsilon}b_{n}^{\beta}$. 
\begin{enumerate}
\item If $\varepsilon=1$ we use two times (C1) and we are done.
\item Otherwise use either (C$\alpha\beta$) 
and then (C1) if $\alpha\neq\beta$ or else both 
(C12) and (C21) and then (C1), if $\alpha=\beta=2$. 
\end{enumerate}
\end{enumerate}
\end{enumerate}
This proves that every vertex descends to $V_n^0$. But these    
vertices have not outgoing edges, as can be easily seen. 
When we use the unoriented edges some new vertices have to be added. 
But it is easy to see that these new vertices do not have    
outgoing edges either. Since any vertex
has a semi-oriented path ending in $V_n^0$ our claim follows. \end{proof}

\begin{remark}  \label{rem:loop}  
The moves (C12) and (C21) are really necessary 
for the conclusion of proposition \ref{lem:3} hold true.    
For instance look at the case $\al=\be=0$.   
From $b_{j+1}b_j^2b_{j+1}^2$ only (C2) can be applied; 
its reduction is a linear combination containing 
the factor $b_{j+1}^2b_j^2b_{j+1}$.  
If we continue, then we shall find at each stage 
one of these two monomials.  Moreover, after making all  
possible reductions at the    
second stage, we recover the word 
$b_{j+1}b_j^2b_{j+1}^2$. Therefore  there exist    
closed oriented loops in the graph. 
In particular the connected component of $b_{j+1}b_j^2b_{j+1}^2$     
has no minimal element, unless we enlarge the graph by adding 
the extra edges  associated to (C12) and (C21).   
For general $\al, \be$ a similar argument holds and it can be checked  
by a computer program.   
If one does not use (C12) or (C21) then the reduction process   
for  $b_{j+1}b_j^2b_{j+1}^2$ yields at the sixth stage    
a  sum of words generating  an oriented loop. \end{remark}

We are able now to define the bicolored graph 
$\Gamma_n^*(H)$, where the non-uniqueness 
of the reduction process  is measured by means 
of an ideal $H\subset R$. 
\begin{definition}
Consider a minimal vertex of $\Gamma_n$ which can
therefore be written as the linear combination:
$v=\sum_{i, k}\lambda_{i,k} (x_{i,k} b_n^{k} y_{i,k})$, 
where $k\in\{0, 1, 2\}$, $x_{i,k},y_{i,k}$ are words from 
$F_{n-1}$ and $\lambda_{i,k}$ are scalars. 
Then we join $v$ by an oriented red edge
to the vertex of $\Gamma_{n-1}$ which corresponds to the 
linear combination: 
\[w=\sum_{i}\lambda_{i,0} (x_{(i,0)} y_{i,0)})+ 
\sum_{i} z \lambda_{i,1} (x_{i,1} y_{i,1}) +
\sum_{i} t \lambda_{i,2} (x_{i,2} y_{i,2})\]    
Finally, the level zero graph  $\Gamma_0(H)$ is the graph having the 
vertices corresponding  to the module $R$. 
Two vertices  of $\Gamma_0(H)$ are connected by an unoriented edge 
if the corresponding elements lie in the same coset of $R/H$, 
where $H$ is a given ideal of $R$. 
\end{definition}
\begin{remark}
The submodule $H$ is necessary because going  
on different descending paths, we might obtain different elements
of $R$. 
\end{remark}   

\section{The coherence conditions for $\Gamma_n^*(H)$}
\subsection{General considerations}
The purpose of this section is to reduce the  
coherence test for $\Gamma_n^*(H)$ 
to finitely many algebraic checks.   
   
We test  the coherence conditions for each 
$\Gamma_n^*(H)$ by recurrence on $n$.    
Notice that for $n\in\{1,2\}$ there are no non-trivial
requirements for $H$.    

The coherence  test for  $\Gamma_n$ (fixed $n$) amounts 
to checking that all open pentagon configurations, which are 
infinitely many, verify (PC). 
Moreover the open pentagon configurations 
themselves can be organized in a pattern which has 
the additional structure of an algebra, in fact a planar 
algebra. We will not make use directly of this 
algebra structure in the sequel. 
However, it can be inferred from it that it is 
enough to verify the (PC) only for those o.p.c. 
which generate this algebra. 
A detailed analysis of these generators 
reduces then the test problem to an explicit 
infinite family of o.p.c. 
At this point we notice that the (PC) might not hold 
for all o.p.c. in this family. 
Now, one enlarges $\Gamma_n$ to the tower of  
colored  graphs $\Gamma^*_n$ and look for the weaker
(CPC) condition for the last one. Eventually, we show that   
the (CPC) for these graphs can be reduced 
to finitely many checks.

The o.p.c.   
$\left[w_0,w_1,\dots,w_{m+1}\right]$ is said to be {\em irreducible} 
if none of the vertices $w_1,w_2,\dots,w_{m}$ has an 
outgoing edge (except the obvious one for $w_1$ and $w_{m}$).    
\begin{lemma} \label{lem:3b}   
\begin{enumerate} 
\item In order to verify (PC) it suffices to restrict to irreducible configurations.    
\item It suffices to verify (PC) only for words from $F_n$.    
\item  Let $\left[w_0,w_1,\dots,w_{m+1}\right]$ be an o.p.c. and 
$w_j'=Aw_jB$, for $j\in\{0,\dots,m+1\}$, where $A,B$ are two arbitrary
words.  
If (PC) holds for $\left[w_0,w_1,\dots,w_{m+1}\right]$, then 
it holds for $\left[w_0',w_1',\dots,w_{m+1}'\right]$.       
\item Suppose that (PC) holds for the two o.p.c. 
$\left[w_0,w_1,\dots,w_{m+1}\right]$ and 
$\left[y_0,y_1,\dots,y_{k+1}\right]$. Then for all $A, B, C$ the    
(PC) is valid also for the following mixed o.p.c.:
   
$\left[Aw_0By_1C,Aw_1By_1C,\dots,Aw_{m}By_1C,Aw_{m}By_2C,\dots, 
Aw_{m}By_{k+1}C\right]$.    

More generally, if one keeps fixed the endpoints of 
the o.p.c., then we can mix the unoriented edges of 
each subjacent o.p.c. following an arbitrary pattern. Specifically,    
let $(i_s,j_s)\in \{0,1,\dots,m+1\}\times\{0,1,\dots,k+1\}, 
s\in\{1,\dots,p\}$ such that:    
$i_0=0<i_1\leq i_2\leq\dots\leq i_p, j_p=k+1>j_{p-1}\geq\dots\geq 0$ 
and  $i_{s+1}-i_s+j_{s+1}-j_s=1$, for all $s$. 
Then the o.p.c.    
$\left[Aw_{i_0}By_{j_0}C,Aw_{i_1}By_{j_1}C,\dots,Aw_{i_p}By_{j_p}C\right]$  fulfills the (PC).    
\end{enumerate}

 \end{lemma}    
\begin{proof}  1) First, any o.p.c.  can be decomposed into irreducible ones. 
Further, if each irreducible component 
satisfies the (PC) then their composition verifies, too.

2) The reduction transformations  acting on different 
monomials of a linear combination commute with each other.
   
3) Obvious.   
   
4) The simplification transformations for 
$w_{m}$ and $y_1$ commute with each other. \end{proof}

\noindent From now on we can restrict  ourselves to 
analyze only those o.p.c. $[w_0,w_1,\dots,w_{m+1}]$ which are 
irreducible. 

\subsection{Resolving the diamonds}

We consider first the case when 
the top line is trivial i.e. $m=1$  and so 
the pentagon degenerates into a diamond.

\begin{lemma} \label{lem:topline}  
If the top line is trivial  then the (PC) holds.   
\end{lemma}    
\begin{proof}   
By using lemma \ref{lem:3b} there are only finitely many words 
$w$ on the top line, to check.  Furthermore $w=abc$, 
where $ab, bc\in \{b_{j+1}^3, b_{j+1}b_jb_{j+1}, 
b_{j+1}b_j^2b_{j+1}, b_{j+1}^2 b_j^2 b_{j+1}, 
b_{j+1} b_j^2 b_{j+1}^2 \}_{j\in\{1,\dots,n-2\}}$.  
The number of cases to study  can be easily reduced, since:  
\begin{enumerate}  
\item If $b$ is the empty word, then the (PC) holds;  
\item By homogeneity it suffices to consider $j=1$; 
\item  Let $w^*= w_r\cdots w_1$ denote the reversed word 
associated to $w= w_1\cdots w_r$.  If the (PC) holds  for $w$, 
then it also holds for $w^*$; 
\item Several cases, as  $b_{j+1}^3 b_j b_{j+1}$, 
can be easily tested at hand.  
\end{enumerate}  
The nontrivial situations are those when a (C12)-move (and then a   
(C2)-move) can be applied.  
It suffices therefore to check the case of  
$b_{j+1}^2 b_j^2 b_{j+1}$, since 
$b_{j+1}b_j^2b_{j+1}^2$ is its reversed and the remaining 
$b_{j+1}^{\varepsilon_1}b_j^2 b_{j+1}^{\varepsilon_2}$   
($\varepsilon_i\in\{2, 3\}$) are consequences of these two.
Then we have the situation depicted in the diagram: 
\[ b_2S_1 \longleftarrow b_2^2b_1^2b_2 \longrightarrow C_1\]
where $S_1, C_1$ are those from (5-6).  
If we apply (C12) and (C21) whenever it is possible on  
$b_{2} S_1$,  then after  
a long computation we find a  common minimal element  
for $b_2S_1$ and  $C_1$.   
\end{proof}

\begin{remark} We used a computer code in order to obtain 
the complete  oriented graph associated  to  the reductions 
of $b_{2}^2 b_1^2 b_{2}$:

\hspace{80pt}\psfig{figure=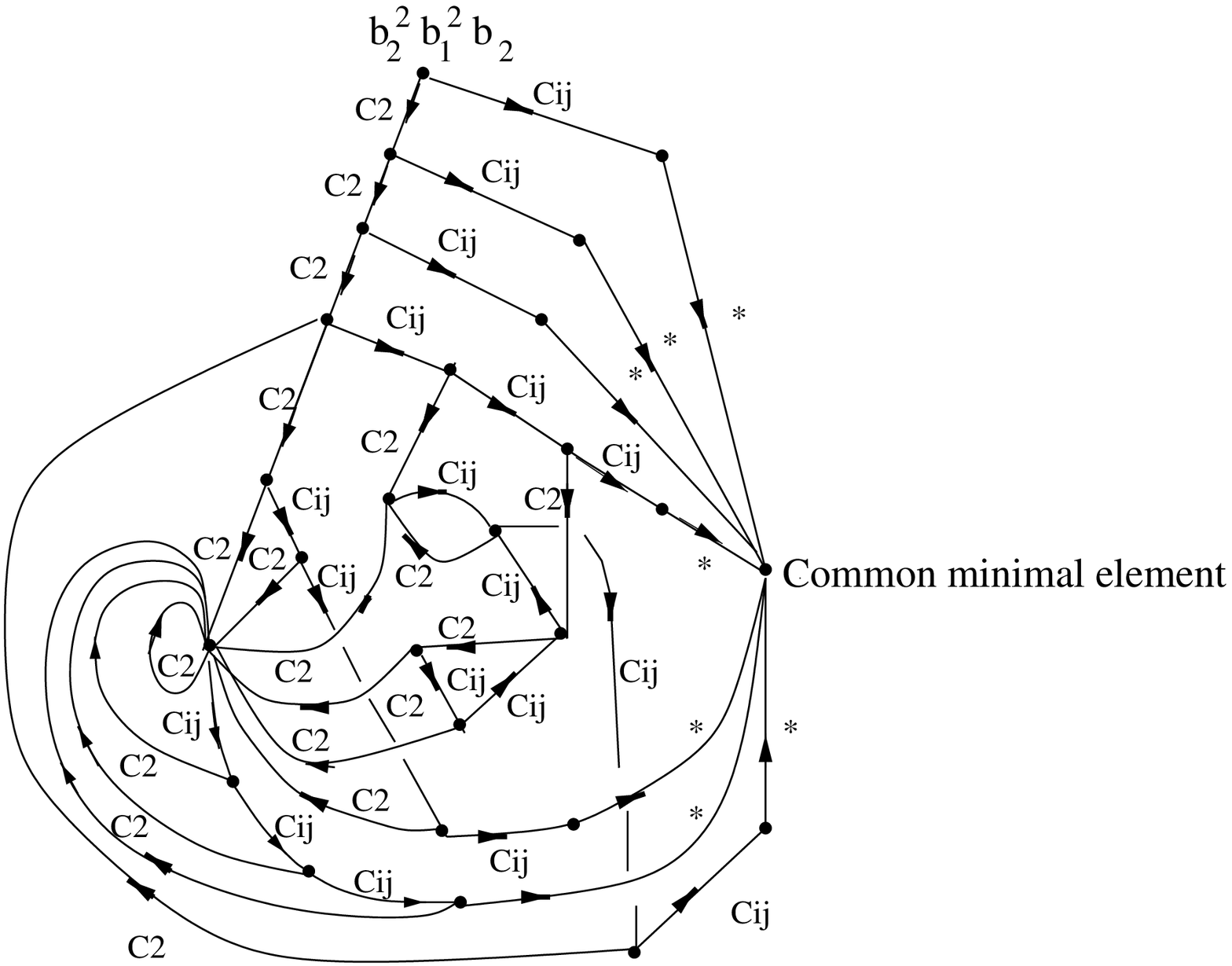,width=11cm} 
    
\noindent Its  vertices are linear combinations in  
words in $b_1$ and  $b_{2}$. The edges are 
labeled by the corresponding reduction. When there are no     
sub-words $b_{2}^2 b_1^2 b_{2}$ or $b_{2} b_1^2 b_{2}^2$
in the factors  of a vertex, its   
reduction is unique;  we marked then the
respective edges by an asterix. The label (Cij) stands for 
the convenient one among (C12) and (C21).  
As we already noticed in Remark \ref{rem:loop}, 
if we apply six times the simplification procedure  
without the use of (Cij)'s then we find a closed loop.  
\end{remark}

\subsection{The diagrams associated to o.p.c.}   
We will be concerned henceforth with the o.p.c.  
having  nontrivial top lines.  
By Lemma \ref{lem:topline} we can suppose that  
$w_1$ and $w_{m}$ have each one exactly one outgoing edge. 
Moreover, an o.p.c. is determined by the following data: 
\begin{enumerate}
\item The word $w=w_1$. Assume that $w$ has length $k$.  
\item The sequence  $w_1,\dots, w_m$, which is 
encoded in a permutation $\sigma\in S_{k}$,
with a specified decomposition into 
transpositions. 
\item The two reduction transformations which simplify 
$w$ and respectively $w_m$.  These should also determine 
uniquely the blocks of letters in $w$ and $w_m$ to which 
the transformations apply.  
\end{enumerate}

Set $T_j$ for the transposition which  interchanges 
the letters on the positions $j$ and $j+1$. 
Let $P(w)$ denote the set of those permutations which 
can be realized on the top line of an o.p.c. having 
its left upper corner labeled $w$. Permutations from $P(w)$ 
will be called permitted permutations. One can  
characterize them as follows. 
Let  $e_w:\{1,2,\dots,l\}\longrightarrow \{1,2,\dots,n-1\}$ 
denote the evaluation map:     
   
\[e_w(j)= \mbox{ the index of the letter lying on the } j\mbox{-th 
position in }  w  \]
Recall that the index of $b_j$ is $j$.     
Consider $\sigma\in P(w)$. Then the permutation 
$T_j\sigma$ is also  permitted if and only if the following 
inequality holds true: 
\[ |e_{\sigma(w)}(j)-e_{\sigma(w)}(j+1)|>1 \]
 
\begin{definition}   
Two permitted permutations $\sigma$ and $\sigma'$ from $P(w)$,
together with their 
specific decomposition into transpositions, 
are  said to be equivalent if the (PC) holds 
true or fails for their associated  o.p.c.,  simultaneously.    
\end{definition}

\begin{lemma} \label{lem:4}   
\begin{enumerate}
\item  Suppose that $\sigma_1T_jT_i\sigma_2\in P(w)$, $|i-j|>1$.    
Then $\sigma_1T_iT_j\sigma_2\in P(w)$ and these two permutations are equivalent.    
   
\item  Suppose that  $\sigma_1T_{i+1}T_iT_{i+1}\sigma_2\in P(w)$. Then    
$\sigma_1T_{i}T_{i+1}T_{i}\sigma_2\in P(w)$ and these two permutations are equivalent. The converse is still true.    
   
\item  If $\sigma_1T_i T_i\sigma_2\in P(w)$ then $\sigma_1\sigma_2$ is permitted and equivalent to    
the previous one.    
\end{enumerate}
\end{lemma}   
\begin{proof} The existence in the first case is equivalent to    
 $| e_{\sigma_2(w)}(j)-e_{\sigma_2(w)}(j+1)|>1$    
 and $| e_{\sigma(w)}(i)-e_{\sigma(w)}(i+1)|>1$,    
 so it is symmetric.    
In the second case also it is equivalent to    
$| e_{\sigma_2(w)}(j+\varepsilon_1)-e_{\sigma_2(w)}(j+\varepsilon_2)|>1$ for all $\varepsilon_j\in\{0,1,2\}$,    
so it is again symmetric. The equivalence is trivial. \end{proof}

\begin{corollary}   
Two different decompositions into transpositions 
of the permutation $\sigma$ lead to equivalent 
o.p.c.
\end{corollary}   
   
We will use  a graphical representation for 
the decomposition of $\sigma$ into    
transpositions, similar to the braid pictures (see picture 
below).  We specify on the top and bottom lines of the rectangle   
the values of the respective evaluation maps. 
Further, the diagram is made of arcs which 
connect the  points on the top to the points 
on the bottom having the same indices; these arcs will be called 
trajectories, or strands in the sequel. We denote by $e(w)$ 
the vector $(e_w(j))_{j=1,\dots,k}$, which can be seen
as a word in the free group (monoid) on $n-1$ letters.   
   
\hspace{158pt}\psfig{figure=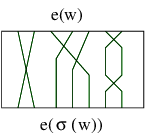,width=3cm}   
   
\noindent 
This picture will be called a {\em diagram} of the 
respective o.p.c.  
Notice that the strands in a diagram inherit 
a labeling  by  the common indices of their endpoints. 
There is also a natural orientation on them, going from the 
top to the bottom.  

The reduction blocks are sets of consecutive 
endpoints of strands (from three to five) in the 
upper and lower lines of a given diagram, corresponding 
to the sub-words on which the 
simplification transformations acts.  
We call them accordingly, the top and the bottom block.

We will draw  below the incomplete diagram
consisting only of  those trajectories of  the six (to ten) 
elements which enter in    
the two reduction  blocks.

\begin{example}
Suppose for instance that the reductions 
consist of two transformations of type (C0). 
This implies that $e(w)=xiiiy$ and $e(\sigma(w))=x'jjjy'$. 
\begin{enumerate}
\item Assume that $i=j$. Then the trajectories    
of the $i'$s can be assumed to be 
disjoint since the transposition which invert 
the letters in the couple $ii$ has trivial 
effect when looking at the word $w$ and its transformations.     
Thus the possible trajectories of these six strands fit 
into the four cases, according to the number of strands
connecting the upper and lower blocks, which might be 
0, 1, 2 or 3.      
\item Further, if $i\neq j$ we have again two sub-cases. 
\begin{enumerate}
\item If $|i-j|=1$  then the trajectories labeled 
$i$  must be disjoint from those labeled $j$, and hence there is only
one obvious combinatorics.   
\item If $|i-j|\neq 1$ then there are sixteen diagrams 
up to isotopy (see \cite{Fun} for a list). 
\end{enumerate}
\end{enumerate}
\end{example}

\begin{remark}  
One can describe all configurations of the strands 
involved in  a pair of reduction transforms (C1)-(C0), 
(C2)-(C0), (C12)-(C0), (C21)-(C0) (see \cite{Fun} 
for an exhaustive list), similar to that from the example above.      
\end{remark}

\begin{definition}  
A diagram is called {\em interactive} if there is at least
one strand connecting the upper and lower blocks. 
\end{definition}

\begin{lemma}    
The  (PC) holds true for the o.p.c. associated 
to non-interactive diagrams.    
\end{lemma}    
\begin{proof}  We call the strands which come or arrive to the reduction blocks 
essential strands. 
\begin{enumerate}
\item 
If the essential arcs coming from the top block are disjoint from 
those arriving in the bottom block then $w=xy$, $\sigma(w)=x y'$, 
where the first block is contained in $x$ and the second one 
in $y'$. These two reductions commute with each other. 
\item If there is an essential strand labeled $i$ which 
intersects some essential strand of the other block, 
then it will intersect  all of them. In particular $b_i$ 
commutes with all letters of the reduction block. 
Moreover, a simple verification shows that, if  
$b_i$ commutes with all letters of the monomial 
from the left hand side of one formula among (3)-(7), 
then it will commute with the elements from the 
right hand side of the same formula. This shows that the 
commutations depicted in the diagram can be realized
after the first reduction transformation (of the upper block).    
This implies our claim. 
\end{enumerate}
\end{proof}

Therefore,  it remains to understand the interactive configurations.

\begin{lemma} 
It suffices to check the (PC) for those 
interactive configurations whose essential strands 
are as following: 

\hspace{92pt}\psfig{figure=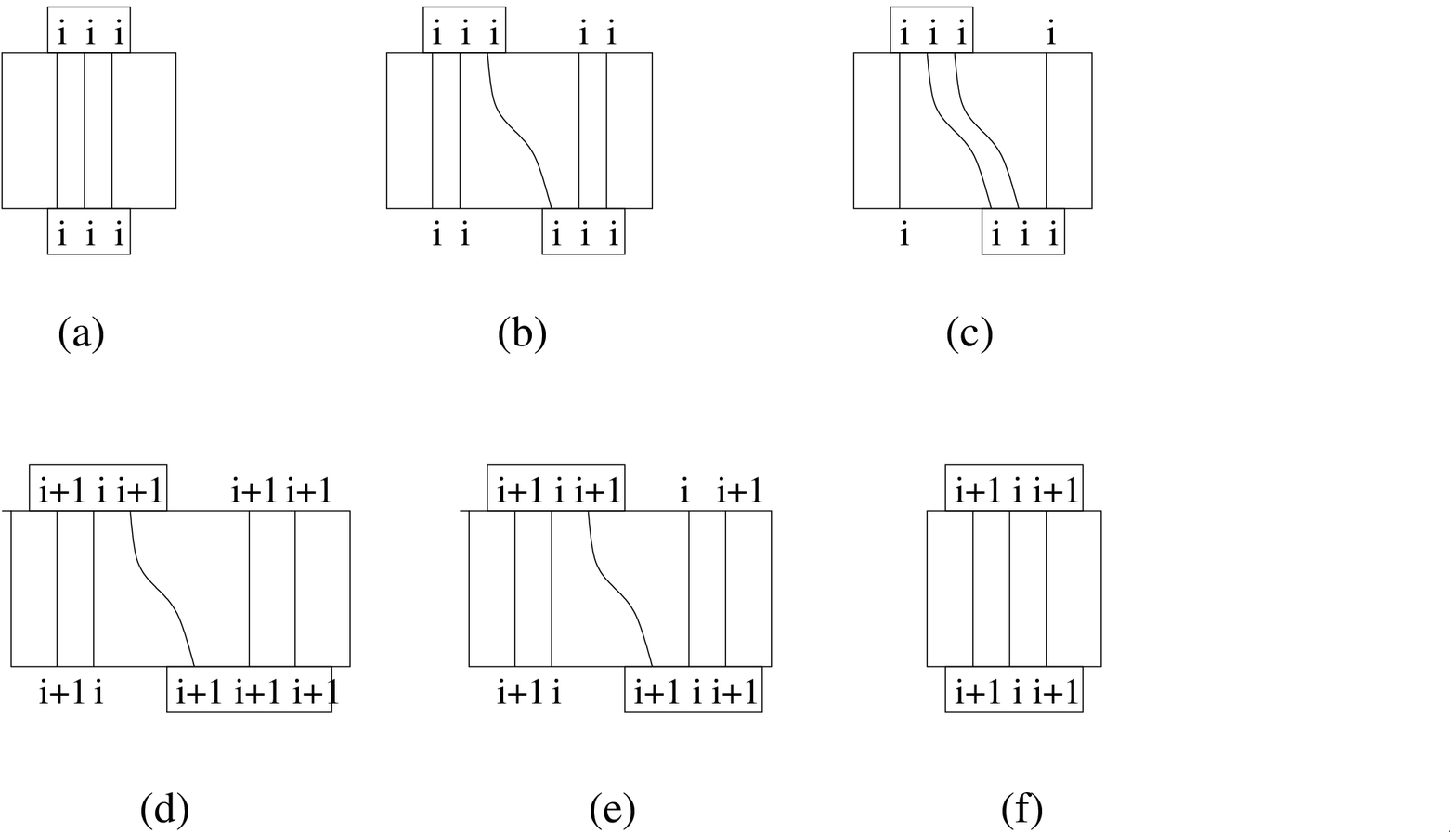,width=10cm}   
\end{lemma}
\begin{remark}   
We represented in the picture above 
each block as a sequence of three letters, but some 
of the letters are allowed to have exponent 2, and thus to represent 
two letters in a genuine diagram. Moreover, in this situation 
we require that the two strands coming from two consecutive letters 
labeled by the same index be parallel, and thus to arrive on 
two consecutive positions on the bottom line. Therefore, 
the couple of 2-strands can be identified with one strand
in the picture above. 
\end{remark}
\begin{proof} There are no restrictions arising 
from the above identification of two
parallel strands because their labels are the same. 
This means that any permutation involving 
one of the two strands  is also  allowable for the second one, 
as well. Thus, we can always get such a normal form for the 
respective interactive configuration. 
\end{proof}

\subsection{Solving the o.p.c. associated to non-interactive diagrams} 
   
The (PC) is verified in the cases (a),(b),(c),(d) and (f) 
by direct computation of the first step of 
their respective simplifications. 
The only relations needed are the consistency of relations 
defining  the algebra $K_3(\al, \be)$. We skip the details.

Let us check a sub-case of (d), corresponding to 
the pair of transformations (C$\varepsilon$)-C(0), 
where $\varepsilon\in\{0,1,2\}$.  
The monomial to be reduced has the form    
$w=b_{i+1}b_i^{\varepsilon}b_{i+1}xb_{i+1}^2$, which is 
weakly equivalent in the o.p.c.   
to $w'=b_{i+1}b_i^{\varepsilon}xb_{i+1}^3$.  We write 
below $\sim$ for the weak equivalence of words. 
Notice that  all 
letters  of $x$ should commute with $b_{i+1}$ because the
respective strands will cross each other.  
Thus we may suppose  that $x$  lies in $F_{i-1}$. 
Therefore:  
$x\rightarrow x_0b_{i-1}^{j_1}b_{i-2}^{j_2}\cdots b_{i-p}^{j_p}$, with
$x_0\in F_{i-2}$. Again, we can restrict ourselves   
to  the situation when $x_0=1$.  
Consider now the case $\varepsilon=2$ because 
the other cases  are trivially verified. 
Set $q=b_{i-2}^{j_2}\cdots b_{i-p}^{j_p}$.  
We have then the following  reduction transformations:   

\[ S_jb_{i-1}^{j_1}b_{i+1}^2q\longleftarrow w \sim w'\longrightarrow
 b_{i+1}b_i^{2}E_jb_{i-1}^{j_1}q\] 
where $S_j, E_j$ as above. From the lemmas 4.1 and 4.2 
it follows that the (PC) holds for:  

\[S_jb_{i+1}^2b_{i-1}^{j_1}q\longleftarrow 
b_{i+1}b_i^{2}b_{i+1}^3b_{i-1}^{j_1}q \longrightarrow
b_{i+1}b_i^{2}E_jb_{i-1}^{j_1}q \]
Since  $S_jb_{i-1}^{j_1}b_{i+1}^2q$ is weakly 
equivalent to $S_jb_{i+1}^2b_{i-1}^{j_1}q$, we are done.  

All remaining cases but (e) follow by  similar computations.    
However, for the diagrams of type  (e) the situation is different. 
Using the commutation rules as above one must preserve 
the term $b_{i-1}^{j_1}$. So, we have to check the configurations
where the word $w$ is given by:

\[w=xb_{i+1}^{\alpha}b_i^{\varepsilon}b_{i+1}^{\beta}b_{i-1}^{\mu}b_{i}^{\delta}b_{i+1}^{\gamma}b_{i-2}^{j_2}\cdots b_{i-p}^{j_p},  
\mbox{ where   } x\in F_{i-1}\] 
At this point one cannot prove  that the (PC) holds for these 
o.p.c. 

\begin{remark}    
In fact the (PC) might not hold since the surjection of 
lemma \ref{lem:man} might have a nontrivial kernel in rank $n=3$. 
\end{remark}

Summarizing what we obtained until now, 
we  proved that these are the only     
o.p.c. that could possibly not verify (PC). 
Moreover, we will check  whether the weaker condition (CPC) 
is valid for these o.p.c.  The explicit computation of 
the minimal elements will show that these are well-defined 
only for the graph $\Gamma_n^*(H)$, for a suitable 
ideal $H$. Let us explain how to find the generators for the 
ideal $H$. 

\begin{proposition}\label{fin}
The (CPC) is verified in $\Gamma_n^*(H)$ if and only if 
it is verified for the following pairs of elements: 
\[  
\, \bt^\xi \, \bd^\epsilon \, \bu^\nu \, \bt^\mu \, \bd^\delta \, \bt^\gamma \; \mbox{and} \;   
 \, \bt^\xi \, \bd^\epsilon \, \bt^\mu \, \bu^\nu \, \bd^\delta \,
 \bt^\gamma \; \; \; \mbox{ for }   
 \xi, \epsilon, \mu,  \nu,  \delta, \gamma \in\{1, 2\}   
\]   
\end{proposition}
\begin{proof}
The only thing one needs to know is that:   
\begin{lemma}    
It suffices to consider  the words $w$ as above with 
$x=1$ and $p=1$.   
\end{lemma}  
\begin{proof}
   
The proposition \ref{indu} shows that any 
admissible functional $\mct$ on $\cup_{n=1}^{\infty}K_{n}(\al, \be)$ 
satisfies:     
   
\[\mct(xuv)=\mct(u)\mct(xv)\; \mbox{ for } x,v\in H(Q,m)
\;\mbox{  and } u\in \langle 1,b_m,b_{m+1},\dots,b_{m+k}\rangle\]    
In the same way one shows that in the simplification process 
the minimal element in $R=\Gamma_0$ 
associated to the word $xuv$ must be the product of 
the minimal elements associated to the  two words $u$ and 
$xv$. 
This proves the claim. 
\end{proof}   
This shows that the cases left unverified can be reduced to 
those which we claimed above.
\end{proof}   

 Thus, the obstructions to the existence of 
the Markov trace come out from these couples.   
 In section 5 we study these obstructions and we find  the    
 ideal $H$ in $R$ containing     
 them.

\section{The computation  of obstructions}  

\subsection{The algorithm}
As we have not yet proved that the trace is well-defined 
we have to specify the choices made 
in the computation of the minimal 
element associated to a given word. 
Moreover, after the verification of the 
(CPC) and commutativity obstructions it will follow 
{\em a posteriori} that all descending paths
in $\Gamma_n^*(H)$ will eventually lead to the 
same element.

Here is the algorithm which was used for computing 
the values of the minimal element in the particular situation 
of the proposition \ref{fin}. Moreover, it can be used 
for any element of the braid group. Notice that the  algorithm 
for reducing the elements of $B_n$ uses recurrently the 
algorithms for the previous stages when simplifying 
elements of $B_{n-1}$. 

\begin{itemize}
\item The input is a word $w$ is an element of $B_n$. 
\item Step 1: use the cubical relations (3) until  
we find a linear combination of words having 
all exponents  within $\{0, 1, 2\}$. We identify each word as an
element of $F_{n-1}$. Further, one writes each word in the form:  
$w=x_1b_{n-1}^{\varepsilon_1}x_2\cdots x_{p}
b_{n-1}^{\varepsilon_p}x_{p+1}$, where $x_i\in F_{n-2}$. 
\item Step 2: if some $x_j$, for $j\in\{ 2,\dots, p\}$ are actually
in $F_{n-2}$ then bring together the two letters 
$b_{n-1}^{\varepsilon_j}$ and $b_{n-1}^{\varepsilon_ {j+1}}$ 
by moving the latter to the left, using the permutations 
(8). 
\item Step 3: perform the steps 1 and 2 
until the output is the same as the output. 
\item Step 4: if $p\geq 2$ then start reducing 
sub-words, starting from left to the right. 
The first sub-word is then 
$b_{n-1}^{\varepsilon_1}x_2b_{n-1}^{\varepsilon_2}$. 
Using the recurrence hypothesis, reduce $x_2$ to a normal 
form in $K_{n-1}(\alpha,\beta)$, and therefore 
write $x_2=y_2 b_{n-2}^{\delta_2} z_2$, where 
$y_2, z_2\in F_{n-3}$.  
Further bring as close as possible the two letters 
$b_{n-1}^{\varepsilon_1}$ and $b_{n-1}{\varepsilon_2}$, by means 
of permutations (8) and obtain the equivalent sub-word 
$y_2b_{n-1}^{\varepsilon_1}b_{n-2}^{\delta_2}
b_{n-1}^{\varepsilon_2}z_2$.
\item Step 5: use the simplification moves  
C(1), C(12), C(2) or C(21), 
according to the values of exponents until we reach an element
where the letter $b_{n-1}$ occurs only once, possibly 
with exponent 2. Consider the new instance of the word 
$w$ by concatenating with the complementary sub-words, 
left untouched. 
\item Step 6: keep repeating the transformations from 
Step 4 until $w$ has a normal form with $p=1$. 
\item Step 7: simplify $w$ by using (9) and 
keep track of the polynomial coefficients. 
If $n=2$ then stop and send as output the coefficients.
Otherwise go to the step 1. 
\end{itemize}
\begin{remark} 
It is important to notice that  the normal form for elements
in $K_3(\alpha,\beta)$ is unique, and hence the step 4 
lead us to a well-defined element for $n=4$. For $n\geq 5$, the 
(CPC) obstructions being verified it follows that the 
output will be independent on the element we chose for the 
normal form at the step 4. 
\end{remark}

\subsection{The CPC obstructions for n=4}   
   
It was pointed out in section 4  that 
the coherence of  $\Gamma_n^*(H)$  depends  
only on the following couples:   
\[
\, \bt^\xi \, \bd^\epsilon \, \bu^\nu \, \bt^\mu \, \bd^\delta \, \bt^\gamma \; \mbox{et} \;   
\, \bt^\xi \, \bd^\epsilon \, \bt^\mu \, \bu^\nu \, \bd^\delta \, \bt^\gamma \; \; \;   
\xi, \epsilon, \mu,  \nu,  \delta, \gamma = 1 \;   
\mbox{or} \; 2   
\]   
Furthermore, if a linear functional $\mct$ is admissible 
then it should verify $\mct(w) = \mct(w^*)$, where 
$w^*$ is the reversal of the word $w$. 
One can therefore reduce ourselves to the study of 
the following 24 couples:     
\begin{itemize}   
\item  $(1.i):$ $\, \bt \, \bd  \, P_i \, \bd^2 \, \bt$ and $\, \bt \, \bd  \, P_i' \, \bd^2 \, \bt  $

\item  $(2.i):$ $\, \bt \, \bd  \, P_i \, \bd \, \bt^2$ and $\, \bt \, \bd  \, P_i' \, \bd \, \bt^2  $

\item  $(3.i):$ $\, \bt \, \bd^2  \, P_i \, \bd \, \bt^2$ and $\, \bt \, \bd^2  \, P_i' \, \bd \, \bt^2  $

\item  $(4.i):$ $\, \bt^2 \, \bd^2  \, P_i \, \bd^2 \, \bt$ and $\, \bt^2 \, \bd^2  \, P_i' \, \bd^2 \, \bt  $

\item  $(5.i):$ $\, \bt^2 \, \bd  \, P_i \, \bd^2 \, \bt^2$ and $\, \bt^2 \, \bd  \, P_i' \, \bd^2 \, \bt^2  $

\item  $(6.i):$ $\, \bt^2 \, \bd^2  \, P_i \, \bd \, \bt$ and $\, \bt^2 \, \bd^2  \, P_i' \, \bd \, \bt $   
\end{itemize}   
\noindent where $P_1= \bu \, \bt, \, P_2 = \bu^2 \, \bt, \, P_3= \bu \, \bt^2 , \, P_4= \bu^2 \, \bt^2,   
P_1'= \bt \, \bu, \, P_2' = \bt \, \bu^2, \, P_3'= \bt^2 \, \bu , \, P_4'= \bt^2 \, \bu^2.$

From now on we denote the difference between the minimal 
elements associated to the left hand side and the right hand side
by the corresponding label $(s,i)$.     
For general  $\alpha, \, \beta$  the computation based on the 
algorithm from above is very  long and    
and we needed to be computer-assisted. 
For more information about the code, see the remark \ref{rem:last}.

One finds  $15$ different    
polynomials from these CPC obstructions, and the 
following identities  among them:   
$(5.2)=-\al (3.2), \;$ $   
(6.2)=\al (1.2), \; $ $   
(1.4)=-\al (1.2)$. Thus,  we must consider the couples   
$(1,2), \,(2,4),$ $ \,(3,2),$ $ \,(3,3), \,$ $(3,4), \,(4,1), \,(4,2), \,   
$ $(4,3),   
 \,(4,4), \,   
(5,3), \,(5,4), \,(6,4)$.   
  
However the 4-variables polynomials we found  above 
should be evaluated at specific values of the 
parameters $(z,\bar{z})$, which are compatible with the 
commutativity requirements for a Markov trace. 
We postpone then the calculation 
of obstructions until  the next section 
where we  find which are the convenient values for 
the parameters, as functions on $(\alpha,\beta)$.

\subsection{Commutativity obstructions}   
We are concerned in this section 
with the commutativity constraints imposed for an admissible 
functional to be a Markov trace:        
\[   
 \mct (a b) = \mct (b a) \; \; \mbox{for all} \; a, \, b  
\]   
\begin{lemma}
An admissible linear functional on $K_3 (\al, \, \be)$ 
satisfies the trace conditions above if and only if 
the the values of $(z,t)$  are given either by the 
type (I) rational parameters:    
\begin{equation}z=\frac{-\be^2 + 2 \al}{\al \be +4}, \;\;\;\; t=    
 \frac{\al^2 + 2 \be}{\al \be +4}\end{equation}  
\noindent or else by the type (II) parameters:     
\begin{equation} 
t = \frac{ 2 \al z - 2 z^2 + \be }{2 + \be z}, \;\;   
\mbox{ where }   z  \mbox{ verifies  }  
 (\al \be + 1) z ^3 + ( \al + \be^2) z^2  + 2 \be z   
 + 1 = 0
\end{equation}
\end{lemma}
\begin{proof}
A  trace $\mct$ defined on $K_3 (\al, \, \be)$
should satisfy the following identities: 
\[
\mct ( \bd \, \bu^2 \, \bd) = \mct ( \bu^2 \, \bd^2) , \;  \mct (   
\bu \, \bd \, \bu^2 \, \bd) = \mct (\bd \,  \bu \, \bd \, \bu^2)   
\]   
These are equivalent to:        
\[   
\mathcal{T}(R_0)= \mathcal{T}(R_1) = 0
\]   
Remark that these are also sufficient conditions for 
an admissible functional 
be actually a trace on $K_3 (\al, \, \be)$. Moreover, the 
equations above can be expressed in the 
following algebraic form:

\begin{eqnarray*}   
0& = &(-\be^3+3 \al \be+4) t^2+(3 \al^2-7 \al \be^2-6 \be+2 \be^4) t+(3 \be^2-\be^5-2 \al-3 \al^2 \be+4 \al \be^3)+ \\
 &   & (2 \al \be^3+\be^2-6 \al^2 \be-10 \al) z t+ 
(-3 \al^3+7 \al^2 \be^2+9 \al \be+4-\be^3-2 \al \be^4) z+ \\
 &  & (3 \al^3\be+7 \al^2-\al^2 \be^3-\al \be^2+2 \be) z^2 \\
0 & = & (\be^2-2 \al) t^2+(4+5 \al \be-2 \be^3) t+ (\be^4-2 \be-3 \al
\be^2+\al^2)+ (2 \be+5 \al^2-2 \al \be^2) z t+ \\
&  & (\be^2+2 \al \be^3-5 \al^2 \be-6 \al) z+ (4+\al^2 \be^2+\al
\be-2 \al^3) z^2\\
\end{eqnarray*} 
The  solutions of these equations are those claimed above. 
\end{proof}  
 
Consider now the following polynomials in $\al$ and $\be$: 

\begin{eqnarray*}
L&=&3\alpha\beta^4+5\alpha^2\beta^5-2\alpha\beta+
  2\alpha^4\beta-7\alpha^3\beta^3-7\alpha^2\beta^2-
  \alpha\beta^7+\alpha^3+   
  (13\alpha^3 \beta^2-10\alpha^2 \beta^4+13\alpha^2    
  \beta-6\alpha\beta^3- \\
& &2\alpha^4+3\alpha+2\alpha   
  \beta^6)t+ (-6\alpha^3 \beta- \al \beta^5   
  -6\alpha^2+3 \alpha \beta^2+5\alpha^2 \beta^3)t^2+   
  (-16\alpha^4 \beta^2-5 \alpha \beta^2-2\alpha^2+ \\
& &3\alpha^5+2 \alpha\beta^5  -13\alpha^3 \beta+11\alpha^3\beta^4-
 2\alpha^2 \beta^6)z+ 
 (-2\alpha\beta^4+15\alpha^4 \beta+2\alpha^2 \beta^5-
  11\alpha^3\beta^3+ 15\alpha^3+6\alpha \beta)zt+\\
& &(-3\alpha-\alpha^3 \beta^5+6\alpha^4 \beta^3-3\alpha^3 \beta^2+
  2\alpha^2 \beta^4- 9\alpha^5   
  \beta-9\alpha^2 \beta-10\alpha^4)z^2
\end{eqnarray*}

\begin{eqnarray*} 
M&=&\alpha-\alpha^4+6\alpha^2\beta -2\alpha^5\beta - 
   2\alpha \beta^3+ 7\alpha^4 \beta^3+  
  11\alpha^3 \beta^2+\alpha \beta^6-7\alpha^2 \beta^4   
 -5\alpha^3 \beta^5+\alpha^2 \beta^7+
  (-21\alpha^3\beta- \\
& &2\alpha^2 \beta^6+2 \alpha\beta^2 +14\alpha^2 \beta^3-
 13\alpha^4 \beta^2-7\alpha^2   
 +10\alpha^3 \beta^4-2 \alpha\beta^5 
 +2\alpha^5)t+(-7\alpha^2 \beta^2   
 +6\alpha^4 \beta+10\alpha^3+\\
& &\alpha \beta^4+\alpha^2   
  \beta^5- 5\alpha^3 \beta^3)t^2+ (-3\alpha^6   
  +2\alpha^3 \beta^6+5\alpha \beta+11\alpha^2   
  \beta^2+16\alpha^5 \beta^2+8\alpha^3+25\alpha^4 \beta   
  -11\alpha^4 \beta^4  \\ 
& &-4\alpha \beta^4-10\alpha^3 \beta^3)z   
 +(11\alpha^4 \beta^3-14\alpha^2 \beta+10\alpha^3   
 \beta^2-\alpha+4\alpha \beta^3-15\alpha^5 \beta  
 -27\alpha^4   
 -2\alpha^3 \beta^5)zt + \\
& &(4\alpha \beta^2-4\alpha^2 \beta^3+\alpha^4 \beta^5+19\alpha^5   
    -\alpha^3 \beta^4+4\alpha^2-3\alpha^4   
    \beta^2 +21\alpha^3 \beta-6\alpha^5 \beta^3
    +9\alpha^6 \beta)z^2
\end{eqnarray*}

\begin{eqnarray*}
N&=&12\alpha^2 \beta^3+\alpha    
 \beta^8-6\alpha^2 \beta^6-2\alpha^2   
 +3 \alpha\beta^2 +11\alpha^3 \beta^4-4   
 \beta^5\alpha-6\alpha^4 \beta^2   
 -7\alpha^3\beta+(-21\alpha^3\beta^3+7\alpha \beta^4+5\alpha^3 \\   
& &+10\alpha^4 \beta-2\alpha \beta^7-2\alpha \beta   
 -17\alpha^2 \beta^2+12\alpha^2 \beta^5)t   
 +(-4\alpha^4+10\alpha^3 \beta^2   
 -3\alpha+\alpha \beta^6+5\alpha^2   
 \beta-6\alpha^2 \beta^4-\\
& &3\alpha \beta^3)t^2+(3\alpha+3\alpha \beta^3+2\alpha^2   
 \beta^7+16\alpha^3 \beta^2-2\alpha   
 \beta^6-7\alpha^4-13\alpha^5 \beta  
 +5\alpha^2 \beta-13\alpha^3 \beta^5   
 +25\alpha^4 \beta^3)z+(\alpha^2-\\  
& &12\alpha^3\beta+10\alpha^5+13\alpha^3   
 \beta^4-\alpha^2 \beta^3-2\alpha^2  
 \beta^6+2 \alpha\beta^5 -24\alpha^4 \beta^2  
 -5 \alpha \beta^2)zt+ \\
& &(5\alpha^3+4\alpha^3 \beta^3+14\alpha^5 \beta^2+8\alpha^4    
\beta+7\alpha^2 \beta^2+\alpha^3 \beta^6   
+5\alpha \beta-2\alpha^2 \beta^5-6\alpha^6   
-7\alpha^4 \beta^4)z^2
\end{eqnarray*}

\begin{proposition}\label{mark}     
Consider an admissible functional $\mct$ 
defined on the whole tower of algebras $\cup_{n=1}^{\infty}K_n(\alpha,\beta)$. 
Suppose that $\mct$ is a trace on $K_3(\al,\be)$. Then 
$\mct$ is a Markov trace on the tower $\cup_{n=1}^{\infty}K_n(\alpha,\beta)$
if the equations: 
\[ L=M=N=0\]
are satisfied. 
\end{proposition}
\begin{proof}
We will prove that the commutativity constraints  
are verified  by induction on $n$.  The claim is true 
for $n=3$, and now we suppose that it holds for all algebras 
$K_m(\alpha,\beta)$, for $m\leq n$.   
In order to prove the claim for $K_{n+1}(\alpha,\beta)$
it suffices to consider  $b\in \{b_1,\dots,b_n\}$ and    
$a$ belonging to some system of generators 
of $K_{n+1}(\alpha, \beta)$, as a module. 
In particular we will choose the set of 
generators $W_{n}$  from section 3.2.    
 
For $b=b_i$, $i<n$ the claim is obvious. 
It remains to check whenever    
$\mct(ab_n)=\mct(b_na)$ holds true.     
There are  three cases  to consider:

\begin{list}{\roman{enumi})}
{\usecounter{enumi}}
\item $a\in K_n(\alpha, \beta)$;
\item $a=xb_ny$, with  $x,y\in K_{n}(\alpha, \beta)$;   
\item $a=xb_n^2y$,  with $x,y\in K_{n}(\alpha, \beta)$.   
\end{list}

\noindent which will be discussed in combination with 
the following six sub-cases:

\begin{enumerate}  
   
\item  $x\in K_{n-1}(\alpha, \beta)$, and $y\in K_{n-1}(\alpha, \beta)$,    
   
\item  $x\in K_{n-1}(\alpha, \beta)$,  and $y=ub_{n-1}v$, $u,v\in K_{n-1}(\alpha, \beta)$,   
      
\item  $x\in K_{n-1}(\alpha, \beta)$,  and $y=ub_{n-1}^2v$, $u,v\in K_{n-1}(\alpha, \beta)$,

\item   $x=rb_{n-1}s$, $r,s\in K_{n-1}(\alpha, \beta)$, $y=ub_{n-1}v$, $u,v\in K_{n-1}(\alpha, \beta)$,     
 
\item   $x=rb_{n-1}s$, $r,s\in K_{n-1}(\alpha, \beta)$, $y=ub_{n-1}^2v$, $u,v\in K_{n-1}(\alpha, \beta)$,   

\item    $x=rb_{n-1}^2s$, $r,s\in K_{n-1}(\alpha, \beta)$,
  $y=ub_{n-1}^2v$, $u,v\in K_{n-1}(\alpha, \beta)$.
\end{enumerate}

\noindent The cases  (*,i), (1,ii) and (1,iii) are trivially verified by 
an immediate calculation.  Furthermore one obtains:

\begin{itemize}
\item[(2,ii)] 
$\begin{array}[t]{lcl}\mct(b_nxb_nub_{n-1}v)& = &tz \mct(xuv)=
\mct(xb_nub_{n-1}vb_n)\end{array}$
\item[(2,iii)] 
$\begin{array}[t]{lcl}\mct(b_nxb_n^2ub_{n-1}v)& = &(\al t+\be z +1) 
\mct(xub_{n-1}v)=(\al t+\be z +1) z \mct(xuv)  \\
&=&\mct(xub_{n-1}b_nb_{n-1}^2v)=
\mct(xb_n^2ub_{n-1}vb_n)
\end{array}$.
\item[(2,iii)]
$\begin{array}[t]{lcl} \mct(b_nxb_n^2ub_{n-1}v)& = &(\al t+\be z +1) \mct(xub_{n-1}v)=(\al t+\be z +1) z \mct(xuv)   
=\mct(xub_{n-1}b_nb_{n-1}^2v)\\
&=&\mct(xb_n^2ub_{n-1}vb_n).\end{array}$.    
\item[(3,ii)] $\begin{array}[t]{lcl} \mct(b_nxb_nub_{n-1}^2v)& = &
t^2 \mct(xuv)=\mct(b_n^2b_{n-1}^2)\mct(xuv)=\mct(b_nb_{n-1}^2b_n)\mct(xuv)=\\   
&=&\mct(xub_nb_{n-1}^2b_nv)=\mct(xb_nub_{n-1}^2vb_n)
\end{array}$
\item[(3,iii)] 
$\begin{array}[t]{lcl}  \mct(b_nxb_n^2ub_{n-1}^2v)& = &(\al t+\be z +1)\mct(xub_{n-1}^2v)=(\al t+\be z +1)t \mct(xuv)\\     
&=&\mct(xuv)\mct(b_{n}^2b_{n-1}^2b_{n})=\mct(xb_n^2ub_{n-1}^2vb_n)
\end{array}$
\end{itemize}

\noindent For the remaining cases, 
we need also to know more precisely the form of $su$.    
Specifically, let us write $su=pb_{n-2}^{\varepsilon}q$ with $p,q\in
K_{n-2}(\al, \be)$ and $\varepsilon \in\{0,1,2\}$.   
We can  show by a direct computation that   
the equalities hold also for (4,ii), (4,iii), (6,ii) and  
(6,iii).   
Moreover,  using  {\sf Maple} we have found that in the 
cases (5, ii) and  (5, iii), for $su = p b_{n-2}^2 q$, 
there are two additional    
identities, which are not consequences of those 
in the previous lemma. The difference 
$\mct(ab)-\mct(ba)$ is expressed as a linear combination 
with polynomial coefficients in 
$\mct (r p b_{n-2}^2 q v), \mct   
(r p  b_{n-2} q v)$ and $ \mct ( r p  q v)$.         
For arbitrary elements $r,p,q,v$ as above the three traces 
above seem to be unrelated. A sufficient condition 
for the commutativity to hold is that the coefficients 
in front of these terms vanish. 
We derive therefore the following obstructions:  
   
\begin{itemize}
\item[(5,ii)]    
$L \mct (r p b_{n-2}^2 q v) + M \mct (r p b_{n-2} q v) +
N \mct (r p  q v)=0$    
\item[(5,iii)] $-\al(L\mct (r p b_{n-2}^2 q v) + M \mct (r p b_{n-2} q v) +
N \mct (r p  q v))=0$
\end{itemize}  
Furthermore the vanishing of $L, M$ and $N$ insures the 
commutativity of the admissible functional. 
\end{proof}

\begin{remark}
It seems that the conditions stated in the proposition \ref{mark}
are also necessary for the existence of a 
Markov trace extension. Nevertheless, the ideal defined by 
all obstructions could not be made strictly smaller, 
even if we could get rid of the equations $L=M=N=0$, because 
of the (CPC) obstructions.  
\end{remark}

\section{The existence of Markov traces}   
   
\subsection{Statements}   
\begin{theorem} \label{first}  
There exists an unique  Markov trace:    
\[   
\mct_{(\al,\, \be)}: \cup_{n=1}^{\infty}K_n(\al, \, \be) \to   
\frac{ \Z [\al, \, \be, \, (4+\alpha\beta)^{-1}]} { ({H_{(\al, \, \be)}})}   
\]   
with type (I) parameters:    
$z=(2\alpha-\be^2)/(\alpha \be + 4)$   
and $\bar{z}=-(\alpha^2+2\be)/(\alpha \be + 4)$,   
where: 
$${H_{(\al\, ,\be)}} = 8\alpha^6- 8\alpha^5\beta^2+2\alpha^4 \beta^4   
+36\alpha^4\beta -34\alpha^3\beta^3+17\alpha^3+8\alpha^2 \beta^5+32\alpha^2\beta^2   
-36\alpha\beta^4 +38\alpha\beta+8\beta^6-17\beta^3+8.$$   
\end{theorem}   
It is more convenient now to   
put $\delta= z^2(\be z +1)$, so that the obstructions 
associated to the type (II) parameters become Laurent   
polynomials in $z$ and $\delta$.     
\begin{theorem} \label{second}  
Set  $\al= -(z^7+\delta^2)/(z^4 \delta)$, $\be=  
(\delta-z^2)/z^3$ and $\bar{z}=-z^2/(\be z+1)=-z^4/\delta$.   
There exists an unique  Markov trace with parameters $(z,\bar{z})$:   
\[   
 \mct^{(z,\, \delta)}:\cup_{n=1}^{\infty} K_{n}(\al, \, \be)\to   
\frac{\Z [z^{\pm 1}, \delta^{\pm 1}]} {(P^{(z,\, \delta)})}   
\]    
where $P^{(z,\, \delta)}=z^{23}+z^{18} \delta-2 z^{16} \delta^2-z^{14} \delta^3-2 z^9 \delta^4+2   
 z^7 \delta^5+\delta^6 z^5+\delta^7$.  
\end{theorem}

\subsection{Proof of Theorem \ref{first}}   
Notice that the parameters  $z, t$ have to satisfy the conditions:   
\[
\mathcal{T}(R_0)= \mathcal{T}(R_1) = 0   
\]   
because the Markov trace vanishes on the ideal $I_3$ of relations
defining $K_3(\al,\be)$. In particular $(z,t)$ are either 
of type (I) or type (II) parameters from (10)-(11).     

For $(z,t)$ as in (10) we derive that $\bar{z} = -t$. 
Set $u:= 1/(\al \be + 4)$, $z_0:=2\al - \be^2$ and $t_0:=\al^2 +  
2\be=:-\bar{z}_0$.

\subsubsection{The commutativity obstructions}   
The equations $L=M=N=0$ are equivalent to: 
   
\begin{itemize}   
\item   $u^2\beta {H_{(\al,\,\be)}}  = 0$    
 \item   $-u^2 (\alpha\beta+2) {H_{(\al,\,\be)}} =0$   
\item  $ u^2(\alpha-\beta^2) {H_{(\al,\,\be)}} =0$   
\end{itemize}   
   
\subsubsection{CPC obstructions}   
   
\begin{itemize}   
\item (1.2): \,  $-u^3\alpha (\alpha-\beta^2){H_{(\al\, , \be)}}W   $   
\item (2.4): \,  $u^3(\alpha-\beta^2) (\alpha^2+\beta){H_{(\al\, , \be)}}W   $   
\item (3.2): \,  $u^3(-\alpha^2 \beta^2+2+\alpha\beta+\alpha^3){H_{(\al, \, \be)}}W   $   
\item (3.3): \,  $u^3(\alpha \beta+2){H_{(\al, \be)}}W $   
\item (3.4): \, $u^3\alpha \beta (\alpha-\beta^2){H_{(\al, \,\be)}}W   $   
\item (4.1): \,  $-u^3(\alpha-\beta^2) (\alpha^2+\beta){H_{(\al,\, \be)}}W   $   
\item (4.2): \,  $u^3\alpha (\alpha^3+2+2 \alpha \beta-\alpha^2 \beta^2-\beta^3){H_{(\al, \, \be)}}W   $   
\item (4.3): \, $u^3\alpha(\alpha^3-\alpha^2\beta^2-2-\beta^3){H_{(\al, \, \be)}}W   $   
\item (4.4): \,  trivial   
\item (5.3): \,  $ -u^3(\beta^2+2 \alpha+2 \alpha^2 \beta){H_{(\al, \, \be)}}W   $   
\item (5.4): \, $u^3\alpha(-\alpha^3 \beta^2-\beta^2-\alpha^2 \beta+\alpha^4){H_{(\al, \, \be)}}W   $   
\item (6.4): \, $-u^3\alpha (\beta+2 \alpha^2) (\alpha-\beta^2){H_{(\al, \, \be)}}W   $   
\end{itemize}   
where  $W=(\alpha+2-\beta)(\alpha^2-2\alpha+4+   
\alpha\beta+2\beta+\beta^2)=\alpha^3+8-\beta^3+6\alpha\beta$.   
   
Consequently, the simplification algorithm defines an admissible 
functional $\mct_{(\al,\be)}$ on the tower of algebras 
$\cup_{n=1}^{\infty}K_n(\al,\be)$  with values in  $\frac{\Z[\al, \, \be, \,
  (4+\alpha\beta)^{-1}]}{(H_{(\al, \, \be)}})$, 
because the (CPC) obstructions vanish. Moreover 
$\mct_{(\al,\be)}$ is a trace on $K_3(\al,\be)$ since the 
parameters verify (10) and it is a Markov trace on the 
whole tower because the commutativity obstruction 
vanish, too.  
  
\subsection{Proof of Theorem \ref{second}}   
Consider now the situation of type (II) parameters. 
We will express all obstructions  as rational   
functions on $z$ and $\be$.

\subsubsection{The commutativity obstructions}   
The equations $L=M=N=0$ are equivalent to: 
\begin{itemize}   
\item   $-Z B_1/(z^7 (z \be+1)^4)= 0$    
\item   $-Z B_2/(z^9 (z \be+1)^5)=0$   
 \item  $Z B_3/(z^7 (z \be+1)^5)=0$   
 \end{itemize}   
  
\subsubsection{The  CPC obstructions}    
\begin{itemize}   
\item (1.2): \,  $-Z B_4 B_5 B_6  /(z^{13} (z \be+1)^8)   $   
   
\item (2.4): \,  $-Z  B_4 B_6  B_7/(z^{15} (z \be+1)^9)   $    
   
\item (3.2): \,  $Z  B_4   B_8/(z^{15} (z \be+1)^9)   $

\item (3.3): \,  $-Z B_4  B_9 /(z^{11} (z \be+1)^7)   $

\item (3.4): \,  $Z B_4 B_5  B_6  \be/(z^{13} (z \be+1)^8)   $

\item (4.1): \,  $ Z  B_4  B_6  B_7/(z^{15} (z \be+1)^9)   $

\item (4.2): \,  $ Z B_4 B_5    B_{10} /(z^{17} (z \be+1)^{10})   $

\item (4.3): \,  $Z B_4  B_5    B_{11}/(z^{17} (z \be+1)^{10})   $

\item (4.4): \, trivial

\item (5.3): \,  $- Z B_4     B_{12}/(z^{13} (z \be+1)^8)   $

\item (5.4): \,  $-Z B_4 B_5    B_{13}/(z^{19} (z \be+1)^{11})   $

\item (6.4): \,  $-Z B_4 B_5  B_6 B_{14} /(z^{17} (z \be+1)^{10}) $

\end{itemize}   
where $Z, B_1, \dots, B_{14}$ are the following polynomials in $z, \be$:   
\begin{enumerate}   
\item $Z= 1+7 z \beta+21  z^2\beta^2+z^3+35 z^3 \beta^3+35 z^4 \beta^4+21 z^5 \beta^5+7 z^6 \beta^6+z^7 \beta^7+z^9 \beta^6+8 z^8 \beta^5+23 z^7 \beta^4+32 z^6 \beta^3+23 z^5 \beta^2+8 z^4 \beta-2 z^6+z^9-z^9 \beta^3-5 z^8 \beta^2-6 z^7 \beta   $   
\item $B_1=3 z^3+z^4 \beta+1+z \beta   $   
\item $B_2=5 z^3+10 z^4 \beta+6 z^5 \beta^2+z^6 \beta^3+4 z^6+2 z^7 \beta+1+3 z \beta+3 z^2\beta^2 +z^3 \beta^3   $   
\item $B_3=\beta+2 z \beta^2+4 z^3 \beta+5 z^4\beta^2 +z^5\beta^3 +z^2 \beta^3 -2 z^5  $   
\item $B_4= (z \beta+ z^2\beta+1+z-z^2) (z \beta+1+2 z^3) (z^4\beta^2 -z^3 \beta^2+z^2\beta^2+1+2 z \beta-z-2 z^2\beta +2 z^2+3 z^3 \beta+z^3+z^4 \beta+z^4)   $   
\item $B_5=1+z^3+z^2\beta^2 +2 z \beta   $   
\item $B_6=z^3 \beta^3+1+2 z \beta+2 z^2\beta^2 +z^3   $   
\item $B_7=1+4 z \beta+6z^2 \beta^2 +2 z^3+4 z^3 \beta^3+z^4 \beta^4+z^6 \beta^3+4 z^5 \beta^2+5 z^4 \beta+z^6)   $   
\item $B_8=z^2\beta^3 +\beta+2 z \beta^2-2 z^2-z^3 \beta   $   
\item $B_9=1+6 z \beta+16z^2 \beta^2 +3 z^3+25 z^3 \beta^3+25 z^4 \beta^4+16 z^5 \beta^5+6 z^6 \beta^6+z^7 \beta^7+3 z^8 \beta^5+13 z^7 \beta^4+24 z^6 \beta^3+24 z^5 \beta^2+13 z^4 \beta+z^7 \beta+z^6+z^9   $   
\item $B_{10}=1+6 z \beta+16z^2 \beta^2 +3 z^3+25 z^3 \beta^3+25 z^4 \beta^4+16 z^5 \beta^5+6 z^6 \beta^6+z^7 \beta^7+z^9 \beta^6+7 z^8 \beta^5+20 z^7 \beta^4+31 z^6 \beta^3+28 z^5 \beta^2+14 z^4 \beta+z^6+z^9+z^9 \beta^3+2 z^8 \beta^2+2 z^7 \beta   $   
\item $B_{11}=6 z \beta+16z^2 \beta^2 +3 z^3+10 z^8 \beta^2+5 z^8 \beta^5+z^7 \beta^7+z^9 \beta^6+12 z^7 \beta+12 z^7 \beta^4+19 z^6 \beta^3+20 z^5 \beta^2+12 z^4 \beta+6 z^6 \beta^6+3 z^9 \beta^3+5 z^6+z^9+1+25 z^3 \beta^3+25 z^4 \beta^4+16 z^5 \beta^5   $   
\item $B_{12}=2 \beta+4z^5 \beta^3 -2 z^5+2z^4 \beta^5 +8 z \beta^2+12z^2 \beta^3 -2 z^2+8 z^3 \beta^4+3 z^4 \beta^2-2 z^3 \beta+z^6\beta^4   $   
\item $B_{13}=1+8 z \beta+29z^2 \beta^2 +63 z^3 \beta^3+80 z^6 \beta^3+29 z^7 \beta^7+13 z^9 \beta^6+17 z^9 \beta^3+91 z^4 \beta^4+57 z^5 \beta^2+23 z^4 \beta+4 z^3+6 z^6+4 z^9+91 z^5 \beta^5+63 z^6 \beta^6+39 z^8 \beta^5+70 z^7 \beta^4+30 z^8 \beta^2+22 z^7 \beta+z^{12}+  
z^9\beta^9 -z^{12} \beta^6+z^{10} \beta^4+2 z^{10} \beta^7+8 z^8 \beta^8-3 z^{11} \beta^5+3 z^{11} \beta^2+7 z^{10} \beta   $   
\item $B_{14}=2+8 z \beta+12 z^2 \beta^2+4 z^3+8 z^3 \beta^3+2 z^4 \beta^4+z^6 \beta^3+6 z^5 \beta^2+9 z^4 \beta+2 z^6   $

\end{enumerate}   
Now the proof follows as above, after    
noticing that $Z(z,\beta)=P^{(z,\, \delta)}(z,\delta)$.

\subsubsection{Corollaries}

\begin{corollary}  
There exist unique  Markov traces:     
\[  
\mct: \cup_{n=1}^{\infty}K_n(0, \, 2\lambda) \to    
\frac{\Z [\lambda]} {(8\lambda^6 - 17 \lambda^3 +1  )}      
\]   
with parameters $z=-\lambda^2$, $t=\lambda$ and $\bar{z}=-\lambda$,
and respectively:   
\[    
\mct: \cup_{n=1}^{\infty}K_n( -2\lambda ,\,0) \to    
\frac{\Z [\lambda]} {(8\lambda^6 - 17 \lambda^3 +1 )}      
\]  
with parameters $z=-\lambda$,  $t=\lambda^2$ and   
$\bar{z}=-\lambda^2$.  
\end{corollary}    
We have a  similar situation   
for the other three solutions.   
In fact for  $\alpha=0$,     
we derive $z=-(t-\be)^2$, where $t$ satisfies $(t^3 - 4\be t^2 +   
5 \be^2 t +1 -2\be^3)=0$. In particular   
$\bar{z}^3-\be \bar{z}^2 +1 =0$ because $\bar{z}=t-\be$.

\begin{corollary}  
There exist unique Markov traces:     
\[    
\mct:\cup_{n=1}^{\infty} K_n\left(0, \, \frac{\lambda^3 +1}{\lambda^2}\right) \to    
\frac{\Z [\lambda^{\pm{1}}]} {(\lambda^9-2\lambda^6+\lambda^3+1 )}  \]  
with parameters   
$z=-\lambda^2$, $\bar{z}=\lambda$ and  
$t=\displaystyle{\frac{2\lambda^3 +1}  
{\lambda^2}}$,   and respectively:   
\[    
\mct: \cup_{n=1}^{\infty}K_n\left(-\frac{\lambda^3+1}{\lambda^2}, \, 0\right) \to    
\frac{\Z [\lambda^{\pm{1}}]} 
{(\lambda^9-2\lambda^6+\lambda^3+1)}\]    
with parameters $z=\lambda$,   
$\bar{z}=-\lambda^2$ and  
$t=\displaystyle{-\frac{2\lambda^3+1}{\lambda^2}}$.   
\end{corollary}


\section{The invariants}    
 \subsection{The definition of $I_{(\al,\,\be)}$ }    
    
As  in section 5.2 we set $z=(2\alpha-\be^2)/(\alpha \be + 4)$,    
 $t=(\alpha^2+2\be)/(\alpha \be + 4)$, $u:= 1/(\al \be + 4)$, $z_0:=2\al - \be^2$ and $t_0:=\al^2 +   
2\be=:-\bar{z}_0$ (notice that in this case $\bar{z} = -t$).

\begin{definition}   
Let $L$ be an oriented link. We set therefore:    
\[ I_{(\al,\,\be)} (L)= \left( \frac{1}{z \bar{z}} \right)^{\frac{n-1}{2}} \left(\frac{\bar{z}}{z    
}\right)^{\frac{e(x)}{2}} \mct_{(\al, \, \be)}(x) \in \frac{\Z [\al, \, \be , \, z_0^{\pm \epsilon/2}, \bar{z}_0^{\pm \epsilon/2}]}    
{(H_{(\al,\,\be)})}  \]   
where $x\in B_n$ is any braid whose closure is  isotopic to    
$L$. Here $\epsilon -1$ is the number of components of $L$ modulo 2.    
\end{definition}   
   
\begin{lemma}  
$ I_{(\al,\,\be)}$ is well-defined.   
\end{lemma}  
  
  \begin{proof}   
Since $b_j^{-1}=b_j^{2}-\al b_j -\be$, we  can suppose that $x$ is a positive braid.  
All the elements in $\Gamma_0(H)$ associated to $x$  are polynomials in  
the variables $z, \, t$   
of degree at most $n-1$.  
The substitutions $z= u  z_0$ and $t= u  t_0$   
imply that,  if $\mct_{(\al, \, \be)}(x)$ and  
$\mct_{(\al, \, \be)}(x)'$ are  representatives of  the trace of $x$, then $\mct_{(\al, \, \be)}(x)'-\mct_{(\al, \, \be)}(x)  
= u^{n-1} G(\al,\, \be) H_{(\al, \, \be)}$,  where $G(\al,\, \be)$   is a polynomial in $\al , \, \be$.  
It follows that:  
$$  
 I_{(\al,\,\be)}(L)= \left(\frac{1}{z_0 \bar{z_0}}\right)^{\frac{n-1}{2}}\left(\frac{\bar{z_0}}{z_0   
}\right)^{\frac{e(x)}{2}} \widetilde{\mct}_{(\al, \, \be)}(x)       
$$   
where we put:
$$\widetilde{\mct}_{(\al, \, \be)}(x)=  u^{-n+1} \mct_{(\al, \, \be)}(x)   
\in \frac{\Z [\al,\,\be ]}{ (H_{(\al, \, \be)})}   $$   
\end{proof}  
   
\subsection{The cubical behaviour}   
   
\begin{proposition} \label{am:ph}    
For any link $K$  there exists some $l \in \{ 0, \, 1, \,2 \}$ such that:  
$$    
I_{(\al,\,\be)}(K) = \frac{\sum_{k \in \Z_+} P_{k}(\be)\al^{k}}{\sum_{k \in \Z_+} Q_{k}(\be)\al^{k}}   
 =\frac{\sum_{k \in \Z_+} M_{k}(\al)\be^{k}}{\sum_{k \in \Z_+} N_{k}(\al)\be^{k}}    
$$    
where $P_{k}, Q_{k}, M_{k}, N_{k}$ are  $(3, k+l)$-polynomials.    
\end{proposition}    
\begin{proof}     
We will show that $M_{k}, N_{k}$ are  $(3, k+l)$-polynomials. The other case is analogous.    
Suppose  first that $x\in B_n^+$, where  $B_n^+$ is the set of positive    
braids  and $n$ is such that  $x\notin B_{n-1}^+$.    
Then $e(x)=|x|$ where $|x|$ denotes the length of $x$.     
In the process of computing the value of the trace on the   
word $x$ we  make two types of    
reductions: either one uses the relations from $K_n(\al, \, \be)$,     
or else one replaces  $a b_l b$ by $ z a b$ (respectively   
$a b_l^2 b$ by $ t a b$). 
In the first alternative the word $y$ is replaced by     
$\sum_{s}(\sum_{k \in \Z_+} D_{k, s}(\al)\be^{k}) y_s$, where 
the $y_s$ are words from  $B_n$,    
the coefficients  $D_{k, s}(\al)$ are $(3,  
k+e(x)-l_s)$-polynomials, and $l_s=|y_s|$.   
In the second case the  word $w$ is replaced by    
$ z w' + t w"$ where $|w'|=|w|-1$ and $ |w"|=|w|-2$.    
When we substitute for $z$ and $t$  their values as functions 
on $\alpha$ and $\beta$   one finds that:    
$$    
\mct_{(\al, \, \be)}(x)= \sum_{k \in \Z_+} u^{s_k}V_{k}(\al)\be^{k}       
$$    
where $0 \le s_k \le n-1$ and $V_{k}(\al) $ are  $(3, k+e(x)  
)$-polynomials. In particular:     
$$    
\widetilde{\mct}_{(\al, \, \be)}(x)= \sum_{k \in \Z_+} u^{s_k-n+1}V_{k}(\al)\be^{k}     
$$    
Now  $u^{s_k-n+1}=\sum_{k \in \Z_+} Y_{k}(\al)\be^{k}$ where $Y_{k}(\al) $ are  $(3, k )$-polynomials.    
Thus, it follows that:     
$$    
\widetilde{\mct}_{(\al, \, \be)}(x)= \sum_{k \in \Z_+} L_{k}(\al)\be^{k}     
$$    
where $L_{k}(\al) $ are  $(3, k+e(x))$-polynomials.

Now, remark that the same reasoning holds true 
for non necessarily positive  $x \in B_n$, by  
getting rid of the negative exponents in $x$ 
by making use of the cubic relation. The only difference 
is that one has to  take into account the 
normalization factor in front of the trace.    
The claim follows.  
\end{proof}

\begin{corollary}    
$I_{(\al,\,0)}(K)=\sum_{i \in \Z_+} a_{3i}\al^{3i}$ and, respectively, $I_{(0,\,\be)}(K)=\sum_{i \in \Z_+} b_{3i}\be^{3i}$,    
where $a_{3i}, \, b_{3i} \in \Z [\frac{1}{2}]$.   
\end{corollary}

\subsection{Chirality and a few other properties of $I_{(\al,\,\be)}$ }

\begin{lemma}    
Set  $x^{\dagger}\in B_n$ for the word  obtained from $x$ by 
replacing  each  term $b_j^{\epsilon}$ by  the corresponding 
$b_j^{-\epsilon}$.  Then the following identity    
${\mct}_{(\al, \,\be)}(x)=    
{\mct}_{(-\be, \,-\al)}(x^{\dagger})    
$ holds true.  In particular, if  the link $K$ is amphicheiral
then  the  identity $I_{(\al,\,\be)}(K) =    
I_{(-\be,\,-\al)}(K)$ is fulfilled.    
\end{lemma}    
\begin{proof}     
Let $Q(b_j)^{\dagger}$ (respectively $R_{0}^{\dagger}$) denote 
the image of  $Q(b_j)$ (and respectively $R_{0}$) 
after  the substitutions $\al \to -\be$,    
$\be \to -\al$ and $b_l \to b_l^{-1}$ for $l=1, \dots, n-1$. It is    
easy to check that $Q(b_j)^{\dagger}=b_j^{-3}Q(b_j) = 0$.    
By some more involved computations 
we verified that  $R_{0}^{\dagger} = R_{1} =0$.   
Since $H_{(\al,\,\be)} =    
H_{(-\be,\,-\al)}$ the claim follows.    
 \end{proof}

\noindent The following properties  have been checked  
by direct calculation (see the table from the appendix).

\begin{enumerate}    
\item $I_{( \al ,\,\be)} $ is independent from HOMFLY and  in  
  particular it distinguishes  among knots that have the 
same HOMFLY    polynomial.   
The knots 5.1 and 10.132   have    
 the same HOMFLY polynomial but different   
$I_{(\al,\,0)}$ and $I_{(0,\,\be)}$ invariants.    
This holds true for the three other couples of prime   
knots with number crossing $\le 10$    
that HOMFLY fails to distinguish, i.e.    
(8.8, 10.129), (8.16,10.156), and  (10.25,\, 10.56).

\item $I_{( \al ,\,\be)} $ detects the chirality of  
those knots with   
crossing  number at most 10,     
where HOMFLY fails i.e. the knots 9.42, 10.48, 10.71, 10.91,   
10.104 and 10.125.  
 
\item The Kauffman polynomial does not detect the  
chirality of  $9.42$ and $10.71$ (see \cite{RGK}). 
Therefore   $I_{( \al ,\,\be)} $ is independent 
   from  the Kauffman polynomial.    
    
\item The $2$-cabling of HOMFLY does not detect the  
chirality of $10.71$ (this result was communicated to us
by H. R. Morton).  Therefore   $I_{( \al ,\,\be)} $ is independent 
from  the $2$-cabling of HOMFLY. 
We notice that the $2$-cabling of Jones polynomial  
can be deduced from Dubrovnic polynomial (\cite{Yam}), which is a variant of Kauffman polynomial 
(\cite{KK}). 
 
\item $I_{( \al ,\,\be)} $ does not distinguish between  
the Conway knot $(C)$ and the Kinoshita-Terasaka knot $(KT)$, 
which form a pair of mutant knots.       
\end{enumerate}    

    
\subsection{The definition of $I^{(z,\, \delta)}$}

\begin{definition}   
For each oriented link $L$ we define:   
\[ I^{(z,\, \delta)} (L)= \left( \frac{1}{z \bar{z}}\right)^{\frac{n-1}{2}}\left(\frac{\bar{z}}{z    
}\right)^{\frac{e(x)}{2}} \mct^{(z,\, \delta)}(x) \,     
 \in \frac{\Z [z^{\pm \epsilon/2}, \delta^{\pm \epsilon/2}]}    
{(P^{(z,\, \delta)})} \]   
where $x\in B_n$ is any braid whose closure is  isotopic to    
$L$ and $\al, \, \be, \, t, \, \bar{z}$ as in Theorem  
\ref{second}. Here $\epsilon -1$ is the number of components modulo 2,  
$\epsilon \in \{1,2\}$.    
\end{definition}

\begin{remark}    
This invariant does not detect the amphichirality  of knots.  
Also $I^{(z,\, \delta)}$ does not distinguish the 
Conway knot from the Kinoshita-Terasaka knot.  
\end{remark}   
  
\begin{proposition}    
The invaraint $I^{(z,\delta)}$ can be expressed as follows:
\[ I^{(z,\, \delta)}(K) = \sum_{k \in \Z} H_{k}(\delta)z^{k}= \sum_{k
   \in \Z} G_{k}(z)\delta^{k} \]  
where $ H_{k}, G_{k}$ are  $(3, k)$-Laurent polynomials.\end{proposition}    
\begin{proof}    
The proof is analogous to the proof of Proposition \ref{am:ph}.   
\end{proof}

\begin{remark}  \label{rem:last}  
 For evaluating obstructions  and traces of  braids we used a 
{\sf Delphi} code.    
The input is a word, or a linear combination of words, 
and we restricted to words representing 5-braids for memory 
reasons.  
One transforms first the word to a sum of  positive words, by using   
the cubic relations. Furthermore 
the transformations $C(j)$ and $C(ij)$   
are used in order to reduce the shape of the word as much as  
possible. When it gets stalked, one allows permutations of the letters.  
The final result is the value of the trace on the braid element.    
The  program is available at:
     
\noindent {http://www-fourier.ujf-grenoble.fr/$\sim$bellinge.html}.    
    
\end{remark}

\vsp

\section{Appendix}   
   
The values of the polynomials for    
$I_{(\al,\,0)}(K)$ and $I_{(0,\,\be)}(K)$    
are displayed below for all knots with no more than 8 crossings.   
The second column is a braid representative for the knot.   
The bold entries in the table are the coefficients of $\alpha^0$  
and, respectively  $\beta^0$. 
The other entries are the non zero coefficients   
of  $\alpha^{3k}$ and    
$\beta^{3k}$ respectively,  for $k \in \Z$.   
For example,   
$$   
I_\alpha(6.2)= [-5 -\frac{19}{4} \al^3 -\frac{1}{2} \al^6],  \;\;  
I_\beta(6.2)= [-16 \be^{-3} +19 -2\be^3]
$$   
The entry ``A" in the last column means that the knot  is  amphicheiral.

$$   
\begin{array}{l|l|l|l|c|}   
\hline   
3.1   & b_1^3  & \mbox{{\bf -1}} \;\; -1/4   &  -8 \;\; \mbox{{\bf 2}} &   \\   
   \hline   
4.1   &   b_1 b_2^{-1} b_1 b_2^{-1}    & 8 \;\;  \mbox{{\bf 10}} \;\;1  &  -8 \;\; \mbox{{\bf 10}}\;\; -1 & A \\   
   \hline   
5.1   & b_1^5  & \mbox{{\bf 0}} \;\; 7/8 \;\;1/8 &  -24\;\; \mbox{{\bf 4}} &  \\    
   \hline   
5.2   &b_1^{2} b_2^{2} b_1^{-1} b_2     & \mbox{{\bf 2}}\;\; 17/8 \;\;1/4  &  -8 \;\; \mbox{{\bf 2}}  & \\   
   \hline   
6.1   &  b_1^{-1} b_2  b_1^{-1} b_3  b_2^{-1} b_3 b_2  &-8 \;\;\mbox{{\bf -16}} \;\;-10 \;\;-1  &  \mbox{{\bf 1}} &  \\   
   \hline       
6.2   & b_1^{-1} b_2  b_1^{-1}b_2^3  &  \mbox{{\bf -5}} \;\;-19/4 \;\; -1/2  & -16\;\; \mbox{{\bf 19}}\;\;-2 &  \\   
   \hline   
6.3   & b_1^{-1} b_2^2  b_1^{-2} b_2  &  \mbox{{\bf -3}} \;\; -1/2  &  \mbox{{\bf -3}}\;\; 1/2 & A \\   
   \hline   
7.1   & b_1^7  &  \mbox{{\bf 0}}\;\;  -5/8 \;\; -9/16\;\; -1/16  & -56 \;\;\mbox{{\bf 8}} &  \\   
    \hline     
7.2   &  b_1^{-1} b_3^3  b_2  b_1^{2} b_3^{-1}b_2 &  \mbox{{\bf -3}} \;\; -11/2\;\; -21/8\;\;-1/4  &  -64\;\;-64\;\;\mbox{{\bf -6}} &  \\   
  \hline   
7.3   & b_1^{2} b_2  b_1^{-1} b_2^4 &  \mbox{{\bf -1}} \;\; -7/4 \;\;-19/16\;\;-1/8  & -64 \;\; 48\;\;\mbox{{\bf -4}} &  \\   
   \hline   
7.4   & b_1^{2} b_2  b_3^2 b_1^{-1}b_2b_3^{-1}b_2&  \mbox{{\bf 0}} \;\; -17/8\;\;-9/4\;\;-1/4   & -64\;\;+128\;\; \mbox{{\bf -78}} \;\;8&  \\   
   \hline   
7.5   & b_1^{4} b_2  b_1^{-1} b_2^{2} &  \mbox{{\bf 0}} \;\; -9/8 \;\;-9/8 \;\;-1/8  &  -24\;\;\mbox{{\bf 4}} &  \\   
   \hline   
7.6   & b_1  b_2^{-1} b_1^{-2} b_3b_2^3b_3&  \mbox{{\bf -4}} \;\;-37/8\;\;-1/2  & -24 \;\;\mbox{{\bf 20}}\;\;-2 &  \\   
   \hline   
7.7   &  b_1 b_3^{-1}b_2b_3^{-1}b_2b_1^{-1}b_2b_3^{-1}b_2& -8\;\;\mbox{{\bf -20}} \;\; -21/2 \;\;-1  &  \mbox{{\bf -19}}\;\;37/2\;\;-2 &  \\   
  \hline   
8.1   & b_1^{-1} b_2 b_3b_2^{-1} b_1^{-1} b_4^{2}b_3b_2  b_4 ^{-1} &   16\;\;\mbox{{\bf 43}} \;\;37\;\;12\;\;1  &  -64\;\;144\;\;\mbox{{\bf -88}}\;\;9 &  \\   
  \hline   
8.2   & b_1^{-1} b_2^5 b_1^{-1} b_2   &  \mbox{{\bf 4}}\;\;59/8\;\;23/8\;\;1/4   &  -24\;\;\mbox{{\bf 36}}\;\;-4 &  \\   
 \hline   
8.3   & b_1^{-2} b_2^{-1}b_1 b_4^{2} b_3  b_4 ^{-1} b_2^{-1}b_3 &  -8\;\; \mbox{{\bf -8}} \;\;-1  &  8\;\;\mbox{{\bf -8}}\;\;1 & A \\   
   \hline   
8.4   &b_1^3 b_3  b_2^{-1} b_3^{-2} b_1b_2^{-1}& 8\;\; \mbox{{\bf 8}} \;\;3/4  & 8\;\; \mbox{{\bf -24}} \;\;  19 \;\; -2&  \\   
   \hline   
8.5   & b_1^3  b_2^{-1} b_1^3  b_2^{-1}&  \mbox{{\bf 1}} \;\; 3 \;\; 19/8 \;\; 1/4  &  -24 \;\;\mbox{{\bf 36}}\;\;-4 &  \\   
   \hline   
8.6  &  b_1^{-1} b_2  b_1^{-1}b_3^{-1}b_2^3 b_3^2  &\mbox{{\bf 5}} \;\; 21/2\;\; 21/4\;\;1/2   &  \mbox{{\bf 1}} &  \\   
   \hline   
8.7   &  b_1^4  b_2^{-2} b_1 b_2^{-1} & \mbox{{\bf 3}} \;\; 9/4 \;\;1/4  &  16 \;\; \mbox{{\bf -25}} \;\; 3&  \\   
 \hline   
8.8   &  b_1^{-1} b_2 b_1^2 b_3^{-1} b_2^2b_3^{-2}& \mbox{{\bf 3}} \;\; 17/4 \;\;1/2  & 16 \;\;\mbox{{\bf -21}} \;\;5/2&  \\   
  \hline   
8.9   & b_1^{-1} b_2b_1^{-3}  b_2^3 &  \mbox{{\bf -7}} \;\; -9 \;\;-1  &  \mbox{{\bf -7}} \;\; 9 \;\;-1  & A\\   
   \hline   
8.10   & b_1^{-1} b_2^2b_1^{-2}  b_2^3& \mbox{{\bf 1}}\;\;2\;\;1/4   &  8 \;\;\mbox{{\bf -8}} \;\;1&  \\   
  \hline   
8.11   &  b_1^{-1} b_2^2  b_3^{-1}b_2 b_3^2 b_1^{-1}b_2& 8\;\;\mbox{{\bf 21}}\;\;147/8\;\;6\;\;1/2   &  -64 \;\; 136\;\;\mbox{{\bf -79}}\;\;8 &  \\   
   \hline   
8.12   &b_1 b_2^{-1}b_3  b_4^{-1}b_3 b_4^{-1} b_2  b_1b_3^{-1}  b_2^{-1} & 24\;\;\mbox{{\bf 44}} \;\;21\;\;2  & -24\;\;\mbox{{\bf 44}} \;\;-21\;\;2 & A \\   
   \hline   
8.13   &  b_1^2 b_2 b_3^{-1} b_2b_1^{-1}b_3^{-2} b_2 & 8\;\;  \mbox{{\bf 12}}\;\;21/4\;\;-1/2   & 8\;\;\mbox{{\bf -28}}\;\;39/2\;\;-2 &  \\   
   \hline   
8.14   & b_1^2b_2^2 b_1^{-1} b_3^{-1} b_2 b_3^{-1} b_2  & \mbox{{\bf 6}}\;\;85/8\;\;21/4\;\;1/2   & -8\;\;\mbox{{\bf 18}}\;\;-2  &  \\    
   \hline   
8.15   & b_1^2  b_2^{-1} b_1 b_3^2 b_2^2 b_3&  \mbox{{\bf 0}}\;\;-17/8\;\;-9/4\;\;-1/4  &  64\;\;-32\;\;\mbox{{\bf 4}} &  \\   
   \hline   
8.16   & b_1^2  b_2^{-1}b_1^2 b_2^{-1}b_1 b_2^{-1}&  \mbox{{\bf -3}} \;\;3/2\;\;1/4  &  \mbox{{\bf -7}}\;\;1 &  \\   
   \hline   
8.17   &  b_1^{-1} b_2  b_1^{-1}b_2^2b_1^{-2}b_2&  \mbox{{\bf -11}} \;\;-19/2\;\;-1  &  \mbox{{\bf -11}}\;\;19/2\;\; -1& A \\   
   \hline   
8.18   &b_1  b_2^{-1} b_1  b_2^{-1}b_1  b_2^{-1}b_1  b_2^{-1}&  -8\;\;\mbox{{\bf -16}}\;\;-10\;\;-1   &  8\;\;\mbox{{\bf -16}} \;\;10\;\;-1& A \\   
   \hline   
8.19   &  b_1  b_2 b_1  b_2b_1  b_2^2 b_1&  \mbox{{\bf 0}}\;\;3/8\;\;1/16  &  64\;\;-64\;\;\mbox{{\bf 1}} & \\   
   \hline   
8.20   &  b_1^3  b_2b_1^{-3}b_2&  \mbox{{\bf 5}}\;\;9/2\;\;1/2  & -8\;\; \mbox{{\bf 0}} &  \\   
   \hline   
8.21   &   b_1  b_2^{-2}b_1^2b_2^3  &\mbox{{\bf 1}}\;\;-1\;\;-1/8   & 8\;\; \mbox{{\bf 0}} &  \\   
   \hline   
\end{array}   
$$   
   

\end{document}